\documentclass[11pt]{article}
\usepackage{type1cm}
\usepackage[english]{babel}
\usepackage[a4paper,margin=1in]{geometry}
\usepackage{amsmath,amssymb,amsthm}
\usepackage{mathtools}
\usepackage{enumitem}
\usepackage{booktabs,array,longtable}
\usepackage[expansion=false]{microtype}
\usepackage{setspace}
\usepackage[skip=0.6\baselineskip plus 2pt]{parskip}
\usepackage{hyperref}
\hypersetup{colorlinks=true,linkcolor=blue,citecolor=blue,urlcolor=blue}

\newtheorem{theorem}{Theorem}[section]
\newtheorem{lemma}[theorem]{Lemma}
\newtheorem{proposition}[theorem]{Proposition}
\newtheorem{corollary}[theorem]{Corollary}
\newtheorem{maintheorem}{Theorem}

\theoremstyle{definition}
\newtheorem{definition}[theorem]{Definition}
\newtheorem{example}[theorem]{Example}
\newtheorem{remark}[theorem]{Remark}

\makeatletter
\renewcommand*\l@section{\@dottedtocline{1}{1.5em}{2.8em}}
\renewcommand*\l@subsection{\@dottedtocline{2}{4.3em}{3.6em}}
\makeatother
\newcolumntype{L}[1]{>{\raggedright\arraybackslash}p{#1}}

\title{Harmonic functions near isolated subanalytic singularities:\\
generalized horns and regular blow-up}
\author{Flavio Mosca}
\date{November 23, 2026}

\begin{document}
\maketitle

\begin{abstract}
\noindent
Let $X\subset\mathbb R^N$ be a subanalytic set of dimension $d\ge2$,
closed near $0$, with an isolated singularity at $0$ and connected link,
whose regular part
carries the induced Riemannian metric. We study harmonic functions on
$X\setminus\{0\}$ near $0$ for the Laplace--Beltrami operator.
\begin{itemize}[leftmargin=2em]
\item If $X$ is bi-Lipschitz equivalent, for the intrinsic metric, to a
\textbf{generalized horn} $\{(z^\gamma y,z):\ y\in L,\ z\ge0\}$, with
$\gamma\ge1$ rational and $L$ a compact connected analytic submanifold of
$\mathbb R^d\setminus\{0\}$, then:
\begin{itemize}[leftmargin=1.5em,label=--,topsep=2pt,itemsep=1pt]
\item \textbf{every bounded harmonic function extends continuously} to $0$;
\item \textbf{every positive harmonic function near {\boldmath$0$} can be
written as {\boldmath$qG_X+u_{\mathrm b}$}}, with $q\ge0$, $u_{\mathrm b}$
bounded and $G_X$ a fixed singular profile;
\item a \textbf{B\^ocher-type theorem} holds: a positive harmonic function
is either bounded or divergent, and two divergent ones have a nontrivial
linear combination that extends continuously.
\end{itemize}
The proof relies on a volume
doubling and Poincar\'e condition on intrinsic balls away from $0$, which
is invariant under bi-Lipschitz homeomorphisms and, together with an
annular connectedness condition, yields the Harnack condition on annuli;
the decomposition follows from the Harnack inequality on the links by a
Picard-principle argument.
\item If $d\ge3$ and the spherical blow-up of $X$ at $0$ is a regular
manifold with boundary in the sense of condition (R), the same
decomposition and B\^ocher-type theorem hold, and the singular profile
$G_X$ has a polyhomogeneous expansion, which is described and may contain
a $\log\rho$ term.
\item For the rotational cusps $|x|^2=z^m$, $m>2$, $d\ge3$, whose
spherical blow-up is not regular, we prove the same decomposition
$qG_X+u_{\mathrm b}$ with a radial profile of order
$\rho^{1-m(d-1)/2}$, and hence the same B\^ocher-type theorem.
The exponent of the profile therefore depends on the rate at which the
sections shrink.
\end{itemize}
\end{abstract}



\clearpage
\tableofcontents
\clearpage

\section{Introduction}
\label{sec:intro}

\subsection{The problem}
\label{subsec:problema}

Let $X\subset\mathbb R^N$ be a subanalytic set of dimension $d\ge2$ with
an isolated singularity at $0$; that is, $0\in X$ and, for some $r>0$,
$X\cap B(0,r)$ is closed in $B(0,r)$ and $X\cap B(0,r)\setminus\{0\}$ is
an analytic submanifold of dimension $d$. Closedness near $0$ makes the
links $X\cap S(0,t)$ compact; it is used throughout
Section~\ref{sec:preliminari}, where the results of Mosca
\cite{Mosca2026} are quoted in this setting, and it cannot be dropped
in Section~\ref{sec:VP}: for the upper half-plane with the origin
adjoined, the positive harmonic function $(x,y)\mapsto y$ violates every
Harnack inequality on balls approaching the omitted boundary. On $X\setminus\{0\}$ we consider the induced Riemannian
metric and the corresponding Laplace--Beltrami operator $\Delta_X$. We
study two properties of harmonic functions near $0$, referred to below as
\emph{removability} and the \emph{B\^ocher property}.
\begin{enumerate}[label=(\alph*)]
\item \emph{Removability.} Every bounded harmonic function on
$X\setminus\{0\}$ near $0$ extends continuously to $0$.
\item \emph{B\^ocher property.} Every positive harmonic function near $0$
is either bounded or tends to $+\infty$ at $0$, and any two positive
harmonic functions tending to $+\infty$ have a nontrivial linear
combination that extends continuously to $0$.
\end{enumerate}
If $0$ is a regular point of $X$, (a) is classical and (b) is B\^ocher's
theorem \cite{Bocher1903}, \cite{AxlerBourdonRamey1992}. If the link has
several connected components, both properties are understood component
by component.

\subsection{Known results}
\label{subsec:noti}

For $d=2$ removability holds for every subanalytic germ with connected
link: a proof based on the Whitney condition and on the modulus of
annuli is given by Mosca \cite{Mosca2026}. It can also be derived from
the quasi-isometry of the germ with a finite union of cones and horns,
proved by Grieser through resolution of singularities
\cite[Theorem~1.2]{Grieser1999}, or from Birbrair's horn theorem
\cite{Birbrair2008}. The B\^ocher property for
$d=2$ follows from the conformal type of the punctured disc
\begin{NoHyper}\cite{Mosca2026}\end{NoHyper}.

For $d\ge3$ the available results concern narrower classes. If $X$ is a
cone over a smooth link, or more generally a stratified space with an
iterated edge metric, H\"older regularity of finite energy solutions is
proved by Akutagawa, Carron and Mazzeo \cite{AkutagawaCarronMazzeo2014};
for metrics quasi-isometric to $|x|^{2\nu}|dx|^2$ on a punctured ball,
H\"older continuity is proved by Wan \cite{Wan2018}.
For cones the B\^ocher property follows from Kondrat'ev's separation of
variables \cite{Kondratev1967}; the case of a spherical link with the
metric $|x|^{2\nu}|dx|^2$, $\nu>-1$, is treated by Jiayu Li and Fangshu
Wan \cite{LiWan2019}. If the germ is intrinsically bi-Lipschitz equivalent
to a ball of $\mathbb R^d$, the decomposition of positive harmonic
functions follows from the theorem of Serrin and Weinberger on
divergence form operators \cite[Theorem~4]{SerrinWeinberger1966}. For
$d\ge3$, none of the regularity results quoted above applies to germs whose sections shrink faster than $\rho$, such as the cusps
$|x|^2=z^m$, $m>2$.

\subsection{Main results}
\label{subsec:risultati}

The first result concerns germs that are bi-Lipschitz equivalent, for the
intrinsic metric, to a generalized horn
\[
X_{\gamma,L}:=\{(z^\gamma y,z)\in\mathbb R^d\times\mathbb R:\ y\in L,\
z\ge0\},
\]
where $\gamma\ge1$ is rational and $L\subset\mathbb R^d\setminus\{0\}$ is
a compact connected analytic submanifold of dimension $d-1$
(Definition~\ref{def:corno-generalizzato}).

\begin{maintheorem}[Theorems~\ref{thm:corni-generalizzati} and~\ref{thm:AD-H}]
\label{mthm:A}
Let $X$ be as above and let $\Psi$ be a homeomorphism from a
neighbourhood of $0$ in $X$ onto a neighbourhood of $0$ in $X_{\gamma,L}$,
with $\Psi(0)=0$, bi-Lipschitz for the intrinsic distances. Then every
bounded harmonic function on $X\setminus\{0\}$ near $0$ extends
continuously to $0$. Moreover $X$ satisfies condition (AD) below, and the
B\^ocher property holds at $0$.
\end{maintheorem}

For $\gamma=1$ the horn is a cone and the result falls within the known
cases. For $\gamma>1$ the proof goes through a condition (VP)
(Definition~\ref{def:VP}) of volume doubling and Poincar\'e type on the
intrinsic balls $B_X(x,s)$ with $s$ small compared with $|x|$; this
condition is invariant under bi-Lipschitz homeomorphisms
(Proposition~\ref{prop:VP-invariante}) and holds for $X_{\gamma,L}$
(Proposition~\ref{prop:corni-VP}). Together with an annular
connectedness condition, (VP) implies the Harnack condition on the annuli
$A(r/b,br)$, from which removability follows. The second assertion follows
from a criterion of independent interest (Theorem~\ref{thm:AD-H}): if the
Harnack inequality holds on the links $L_{r_k}$ along a sequence of radii
$r_k\to0$, with a constant independent of $k$, then the positive harmonic
functions near $0$ that vanish on an outer link form a half-line (Picard
principle), and every positive harmonic function near $0$ is the sum of a
multiple of a fixed singular function and of a bounded harmonic
function. The method is that of Pinchover for Fuchsian type operators
\cite{Pinchover1994}.

The other two results describe the singular profile for $d\ge3$. We say
that $X$ satisfies the \emph{asymptotic decomposition} condition (AD) if
there exists a positive harmonic function $G_X$, tending to $+\infty$ at
$0$ and called the \emph{singular profile}, such that every positive harmonic function near $0$ is of the form
$qG_X+u_{\mathrm b}$ with $q\ge0$ and $u_{\mathrm b}$ bounded
on some punctured neighbourhood of $0$ (Definition~\ref{def:DA}). Under
(AD) and removability the B\^ocher property holds
(Proposition~\ref{prop:picard}).

\begin{maintheorem}[Theorem~\ref{thm:DA-R}]
\label{mthm:C}
Let $d\ge3$. If the spherical blow-up of $X$ at $0$ is a regular manifold
with boundary in the sense of Definition~\ref{def:R} (condition (R)), then $X$ satisfies
(AD) and removability, and hence the B\^ocher property. The profile has
the expansion
\[
G_X=\rho^{-(d-2)}+\sum_{0<k<\ell(d-2)}\mu^{k-\ell(d-2)}c_k(y)
+c_{\log}\log\rho+O(1),
\]
where $\ell\ge1$ is the integer of Definition~\ref{def:R},
$\rho=\mu^\ell$, $y$ is the variable on the link of the tangent cone, the
$c_k$ are smooth functions and $c_{\log}\in\mathbb R$.
\end{maintheorem}

\begin{maintheorem}[Proposition~\ref{prop:cuspide-rotazione}]
\label{mthm:D}
Let $d\ge3$ and let $X=\{(x,z)\in\mathbb R^d\times\mathbb R:\ |x|^2=z^m,\
z\ge0\}$ with $m>2$ rational. Then $X$ satisfies (AD) with a radial
profile $G_X\sim z^{1-\beta}/(\beta-1)$, $\beta:=m(d-1)/2$, and the
B\^ocher property holds.
\end{maintheorem}

Condition (R) holds for cones, for nondegenerate perturbations of
homogeneous cones, for Puiseux-type graphs $\{z=|x|^{a/\ell}\}$ and for complex
germs with smooth projective tangent cone (Remark~\ref{rem:esempi-R}); the
class of Theorem~\ref{mthm:C} is disjoint from that of
Theorem~\ref{mthm:A} with $\gamma>1$. Theorem~\ref{mthm:C} shows that the
singular profile is in general not $\rho^{-(d-2)}$, even when the tangent
cone is smooth away from the vertex: Example~\ref{ex:log} contains a
$\log\rho$ term. Theorem~\ref{mthm:D} shows that the exponent of the
profile depends on the rate at which the sections shrink. In both
theorems condition (AD) also follows from Theorem~\ref{thm:AD-H}; their
content is the description of the profile, on which
Theorem~\ref{thm:AD-H} gives no information.

\subsection{Known results and contributions}
\label{subsec:contributi}

The capacity estimate for annuli, the finite energy property and the
maximum principle of Section~\ref{sec:capacita} are adaptations of
classical results (Remark~\ref{rem:fonti-capacita}). The passage from a
Harnack inequality on annuli to the existence of the limit
(Section~\ref{sec:harnack}) is Moser's argument \cite{Moser1961}.
Continuity for cones falls within the results of Akutagawa, Carron and
Mazzeo (Remark~\ref{rem:ACM}). The Harnack inequality on balls under
volume doubling and Poincar\'e is due to Grigor'yan \cite{Grigoryan1991}
and Saloff-Coste \cite{SaloffCoste1992}, in the local form of Lierl and
Saloff-Coste \cite{LierlSaloffCoste2012}. Volume doubling and Poincar\'e
on remote balls, relatively connected annuli and the resulting Harnack
inequality on annuli are due to Grigor'yan and Saloff-Coste
\cite{GrigoryanSaloffCoste2005} near infinity.

The contributions of this paper are the following.
\begin{enumerate}[label=(\roman*),leftmargin=2.5em]
\item The adaptation to an isolated singular point of the scheme of
Grigor'yan and Saloff-Coste: condition (VP) of Definition~\ref{def:VP},
volume doubling and a weak Poincar\'e inequality on the intrinsic balls
$B_X(x,s)$ with $N_Ps\le|x|/2$, plays the role of their conditions on
remote balls, with $\{0\}$ in place of the region at infinity. We record
its invariance under bi-Lipschitz homeomorphisms of the intrinsic
metric, with explicit dependence of the constants
(Proposition~\ref{prop:VP-invariante}), and the Harnack condition on
annuli under (VP) and (RA$_b$), with a constant depending only on those
of (VP) and (RA$_b$) (Proposition~\ref{prop:VP-H}).
\item The verification of (VP) and (RA$_2$) for the generalized horns
$X_{\gamma,L}$ with $\gamma>1$ (Proposition~\ref{prop:corni-VP}), which
are not covered by the regularity results quoted in
Section~\ref{subsec:noti}, and the resulting Theorem~\ref{mthm:A}.
\item A sufficient criterion for the Picard principle and condition (AD):
the Harnack inequality on the links along a sequence of radii
(Theorem~\ref{thm:AD-H}). The scheme of the proof is that of Pinchover
\cite[Theorem~7.1, Corollary~7.3 and Lemma~8.4]{Pinchover1994}. Once
the Dirichlet problem on the compact annuli, interior elliptic estimates
and the maximum principle at the vertex (Lemma~\ref{lemma:massimo}) are
available, the criterion uses only the Harnack inequality on the links. Together with (ii), it gives the
B\^ocher property for the whole class of Theorem~\ref{mthm:A}.
\item The expansion of the profile for germs with regular spherical
blow-up, including the logarithmic term (Theorem~\ref{mthm:C}).
\item The radial profile of the rotational cusps, whose exponent depends
on $m$ (Theorem~\ref{mthm:D}).
\end{enumerate}
All continuity results, including those for cones, are obtained through
condition (VP).

\subsection{Notation}
\label{subsec:notazioni}

We fix $X$ as in Section~\ref{subsec:problema} and use the following
notation.

{\small
\begin{longtable}{@{}L{3.6cm}L{11cm}@{}}
\toprule
Symbol & Meaning \\
\midrule
\endfirsthead
\toprule
Symbol & Meaning \\
\midrule
\endhead
$B(0,r)$, $S(0,r)$ & open ball and sphere in $\mathbb R^N$ with centre $0$ and radius $r$\\
$\rho(x):=|x|$ & Euclidean distance from $0$\\
$C_0=2$, $C_0^Y$ & constant in the comparison $|x|\le d_X(x,0)\le C_0|x|$ (Lemma~\ref{lemma:curve-radiali}), for $X$ and for a set $Y$ with the same properties\\
$r_0$ & radius fixed in Section~\ref{sec:preliminari}, below which transversality and the product structure hold\\
$X^*:=X\cap B(0,r_0)\setminus\{0\}$ & the smooth part near $0$\\
$L_r:=X\cap S(0,r)$, $L:=L_{r_0}$ & the link at radius $r$\\
$\Omega_r$, $A(s,t)$ & $X\cap\{0<\rho<r\}$ and the annulus $X\cap\{s<\rho<t\}$\\
$g$, $\tau$, $\nabla_X$, $\Delta_X$ & induced metric on $X^*$, Riemannian measure, gradient, Laplace--Beltrami operator\\
$d_X$, $B_X(x,s)$ & intrinsic distance of $X$ and intrinsic ball\\
$\mathcal H^k$ & $k$-dimensional Hausdorff measure\\
$V(r):=\mathcal H^{d-1}(L_r)$ & volume of the link\\
$\mathcal P$ & positive harmonic functions on punctured neighbourhoods of $0$\\
(H$_b$), (H), $C_H(r)$ & Harnack condition on the annuli $A(r/b,br)$, the case $b=4$, optimal constant\\
(VP), (RA$_b$) & volume doubling and Poincar\'e, connectedness of annuli (Section~\ref{sec:VP})\\
(AD), (R), $G_X$ & asymptotic decomposition, regular blow-up, singular profile (Section~\ref{sec:picard})\\
$X_{\gamma,L}$, $h_L$ & generalized horn, induced metric on $L$\\
$\alpha_j^\pm$, $\lambda_j$ & indicial exponents and eigenvalues of the Laplacian of the link\\
\bottomrule
\end{longtable}
}

Constants carry a subscript indicating the condition they come from
($C_H$, $C_V$, $a_P$, $N_P$, $\Lambda_A$, $C_E$); the radii below which
the conditions hold carry the same subscript ($r_H$, $r_V$, $r_A$). The
symbol $C$ without subscript denotes a positive constant whose value may
change from line to line within a proof. Harmonic functions are denoted
by $u$; in Harnack inequalities the nonnegative function is denoted by
$v$.

\section{Preliminaries}
\label{sec:preliminari}

This section collects the facts on the geometry of $X$ near $0$ and the
consequences of the vanishing capacity of the singular point. Unless
stated otherwise, $X$ is as in Section~\ref{subsec:problema}, with
$d\ge2$.

Let $r>0$ be such that $X\cap B(0,r)\setminus\{0\}$ is an analytic
submanifold of dimension $d$. On it we consider the Riemannian metric $g$
induced by $\mathbb R^N$, the Riemannian measure $\tau$, the gradient
$\nabla_X$ and the Laplace--Beltrami operator
\[
\Delta_Xu=\frac1{\sqrt{\det g}}\sum_{i,j}\partial_{y_i}\!\Big(\sqrt{\det g}\,
g^{ij}\,\partial_{y_j}u\Big)
\]
in local coordinates $(y_1,\dots,y_d)$. A function $u\in C^2$ is
\emph{harmonic} if $\Delta_Xu=0$. For $0<r<r_0$, with $r_0$ fixed below,
we set $L_r:=X\cap S(0,r)$, $X^*:=X\cap B(0,r_0)\setminus\{0\}$,
$\Omega_r:=X\cap\{0<\rho<r\}$ and $A(s,t):=X\cap\{s<\rho<t\}$ for
$0<s<t\le r_0$.

\subsection{Uniform transversality and product structure}
\label{subsec:trasversalita-automatica}

For a $C^1$ submanifold $Q\subset\mathbb R^N$ with
$0\in\overline Q\setminus Q$ we define, for $x\in Q$, the Whitney
function $\varpi(x):=\sin\angle(T_xQ,\,\mathbb Rx)$, where
$\angle(T_xQ,\mathbb Rx)$ is the angle between the line $\mathbb Rx$ and
its orthogonal projection on $T_xQ$. The pair $(Q,\{0\})$ satisfies the
Whitney condition at $0$ if $\varpi(x)\to0$ as $x\to0$ in $Q$.

\begin{theorem}[Hironaka \cite{Hironaka1977}]
\label{thm:hironaka}
Let $Q$ be a connected analytic submanifold of $\mathbb R^N$, subanalytic
in $\mathbb R^N$, and let $Z\subset\overline Q\setminus Q$ be a locally
closed subanalytic set. There exists a closed subanalytic subset
$Z'\subset Z$, with $Z\setminus Z'$ dense in $Z$, such that $(Q,Z)$
satisfies the Whitney condition at every point of $Z\setminus Z'$.
\end{theorem}

We apply it with $Z=\{0\}$: since $Z\setminus Z'$ must be dense in $Z$,
we obtain $Z'=\emptyset$, and the Whitney condition holds at $0$,
independently of $d$ and of the codimension.

\begin{lemma}[{Mosca \cite[Addendum~4.7]{Mosca2026}}]
\label{lemma:trasversalita}
If $(X\cap B(0,r)\setminus\{0\},\{0\})$ satisfies the Whitney condition
at $0$, there exists $r_0>0$ such that
\[
|\nabla_X\rho(x)|=\sqrt{1-\varpi(x)^2}\ \ge\ \tfrac12\qquad(x\in X\cap B(0,r_0)\setminus\{0\}).
\]
\end{lemma}

Since $\nabla_X\rho$ is the projection of $\nabla\rho=x/|x|$ on $T_xX$, we
also have $|\nabla_X\rho|\le1$. In what follows $r_0$ is the radius of
Lemma~\ref{lemma:trasversalita}, possibly reduced so that the product
structure and Lemma~\ref{lemma:volume} also hold.
\label{subsec:prodotto}

By the Whitney condition, the restriction $\rho|_{X^*}$ has no critical
points near $0$ (Mosca \cite[Lemma~4.5]{Mosca2026}), and there exists a
diffeomorphism (Mosca \begin{NoHyper}\cite[Corollary~4.6]{Mosca2026}\end{NoHyper})
\[
\Phi:X\cap B(0,r_0)\setminus\{0\}\ \xrightarrow{\ \cong\ }\
L_{r_0}\times(0,r_0),\qquad \Phi(L_r)=L_{r_0}\times\{r\}.
\]
In particular all links $L_r$ are diffeomorphic to $L:=L_{r_0}$. A
diffeomorphism with these properties is obtained, for instance, from the
flow of the vector field $\nabla_X\rho/|\nabla_X\rho|^2$, defined on $X^*$
near $0$, which maps each level set of $\rho$ onto another level set.

\begin{lemma}[{Mosca \cite[Lemma~7.5]{Mosca2026}}]
\label{lemma:volume}
There exists $C>0$ such that $V(r)=\mathcal H^{d-1}(X\cap S(0,r))\le
Cr^{d-1}$ for $0<r\le r_0$.
\end{lemma}

Lemma~\ref{lemma:volume} requires no assumption on $d$. From now on we
assume that the link $L$ is connected; then, by the product structure,
$\Omega_r$ and the annuli $A(s,t)$ are connected.

\subsection{Intrinsic distance from the origin and radial curves}
\label{subsec:distanza-origine}

We denote by $d_X$ the intrinsic distance of $X$, that is, the infimum of
the lengths of rectifiable curves in $X$, and by $B_X(x,s)$ the open ball
with centre $x$ and radius $s$ for $d_X$.

\begin{lemma}
\label{lemma:curve-radiali}
For $x\in X^*$ and $0<t<r_0$ there exists a curve in $X^*\cup\{0\}$
joining $x$ to a point of $L_t$, along which $\rho$ is monotone, of
length at most $2\bigl|\,|x|-t\,\bigr|$; as $t\to0$ it joins $x$ to $0$. In
particular
\begin{equation}
\label{eq:distanza-origine}
|x|\le d_X(x,0)\le C_0|x|\qquad\text{for }x\in X^*,\qquad C_0:=2 .
\end{equation}
\end{lemma}

\begin{proof}
The vector field $\nabla_X\rho/|\nabla_X\rho|^2$ is smooth on $X^*$, and
$\rho$ increases with unit speed along its integral curves, whose speed
is $|\nabla_X\rho|^{-1}\le2$ by Lemma~\ref{lemma:trasversalita}. The
integral curve through $x$, followed in the direction of increasing
$\rho$ if $t>|x|$ and of decreasing $\rho$ if $t<|x|$, reaches $L_t$
after time $\bigl|\,|x|-t\,\bigr|$ and stays in $X^*$, because $\rho$ remains
between $|x|$ and $t$. The first inequality in~\eqref{eq:distanza-origine}
is obvious.
\end{proof}

The constant $C_0$ is denoted separately because it enters the constants
of Section~\ref{sec:VP}. For a set $Y$ with the same
properties as $X$ we write $C_0^Y$; one always has $C_0^Y\le2$.

\subsection{Capacity, finite energy and the maximum principle}
\label{sec:capacita}

Let $\Omega\subset X^*$ be open. For disjoint closed sets
$E_0,E_1\subset\overline\Omega$, the \emph{capacity} of the condenser
$(E_0,E_1;\Omega)$ is
\[
\mathrm{cap}_2(E_0,E_1;\Omega):=\inf\Big\{\int_\Omega|\nabla_Xu|^2\,d\tau:\
u\ \text{Lipschitz on }\overline\Omega,\ u=1\ \text{on }E_0,\ u=0\ \text{on }E_1\Big\},
\]
where $\overline\Omega$ denotes the closure of $\Omega$ in $X^*$; in the
applications $E_0$ and $E_1$ are links $L_{r_1}$, $L_{r_2}$, which are
compact because $X$ is closed near $0$.

\begin{remark}[Sources]
\label{rem:fonti-capacita}
Proposition~\ref{prop:stima-capacita} and Lemmas~\ref{lemma:taglio},
\ref{lemma:energia} and~\ref{lemma:massimo} are not new: they are
adaptations to $X^*$ of classical results, proved here in a
self-contained form. The capacity estimate through the integral of $1/V$
is the upper bound for capacity used in volume criteria for
parabolicity (Grigor'yan \cite[Section~7.1]{Grigoryan1999}).
Removability of sets of zero capacity for bounded solutions goes back to
Serrin \cite{Serrin1964b} in $\mathbb R^n$; for metric measure spaces see
A.~Bj\"orn \cite{Bjorn2006}. In the subanalytic case a related statement,
the density in Sobolev spaces of smooth functions vanishing near the
singular set, is proved by Valette under his hypotheses
\cite[Theorem~2.7 and Corollary~2.8]{Valette2024}; the proofs below do not
rely on it. The only ingredient
specific to subanalytic sets is the estimate $V(r)\le Cr^{d-1}$ of
Lemma~\ref{lemma:volume}.
\end{remark}

For $0<r_1<r_2\le r_0$ let $\Omega:=X\cap\{r_1<\rho<r_2\}$. By Lemma~4.5
of Mosca \cite{Mosca2026} the function $\rho|_{X^*}$ has no critical
points in $\overline\Omega$; the level sets $L_r$, $r\in[r_1,r_2]$, are
therefore compact smooth hypersurfaces of $X^*$, and $V(r)$ is a
continuous positive function on $[r_1,r_2]$.

\begin{proposition}
\label{prop:stima-capacita}
For $0<r_1<r_2\le r_0$,
\[
\mathrm{cap}_2\big(L_{r_1},L_{r_2};\Omega\big)\ \le\
\Big(\int_{r_1}^{r_2}\frac{dr}{V(r)}\Big)^{-1}.
\]
\end{proposition}

\begin{proof}
Let $I:=\int_{r_1}^{r_2}V(t)^{-1}dt$ and
$F(r):=I^{-1}\int_r^{r_2}V(t)^{-1}dt$, so that $F(r_1)=1$, $F(r_2)=0$ and
$F'=-1/(IV)$ is continuous on $[r_1,r_2]$. Since $\rho$ is $1$-Lipschitz
for the intrinsic distance, $u:=F\circ\rho$ is Lipschitz on $\Omega$,
equals $1$ on $L_{r_1}$ and $0$ on $L_{r_2}$, and
$|\nabla_Xu|=|F'(\rho)|\,|\nabla_X\rho|$. By the coarea formula with
respect to $\rho$ and the inequality $|\nabla_X\rho|\le1$,
\[
\int_\Omega|\nabla_Xu|^2\,d\tau=\int_{r_1}^{r_2}F'(r)^2\int_{L_r}
|\nabla_X\rho|\,d\mathcal H^{d-1}\,dr\le\int_{r_1}^{r_2}F'(r)^2V(r)\,dr
=\frac1{I^2}\int_{r_1}^{r_2}\frac{dr}{V(r)}=\frac1I .\qedhere
\]
\end{proof}

\label{subsec:energia}

\begin{lemma}
\label{lemma:taglio}
There exist Lipschitz functions $\eta_k:X\cap B(0,r_0)\to[0,1]$,
$k\ge1$, and radii $0<s_k<t_k$ with $t_k\to0$, such that $\eta_k=1$ on
$\overline{\Omega_{s_k}}$, $\eta_k=0$ outside $\Omega_{t_k}$ and
$\int|\nabla_X\eta_k|^2d\tau\le1/k$.
\end{lemma}

\begin{proof}
Let $t_k:=2^{-k}r_0$. By Lemma~\ref{lemma:volume},
$\int_0^{t_k}dr/V(r)=+\infty$, so there exists $s_k\in(0,t_k)$ with
$\int_{s_k}^{t_k}dr/V\ge k$. Let $F_k$ be the function of the proof of
Proposition~\ref{prop:stima-capacita} on the interval $[s_k,t_k]$; set
$\eta_k:=F_k\circ\rho$ on $A(s_k,t_k)$, $\eta_k:=1$ on
$\overline{\Omega_{s_k}}$ and $\eta_k:=0$ elsewhere. By
Proposition~\ref{prop:stima-capacita} the energy of $\eta_k$ is at most
$1/k$.
\end{proof}

\begin{lemma}[Finite energy]
\label{lemma:energia}
Let $u$ be harmonic and bounded on $\Omega_{r_2}$, with $r_2\le r_0$. For
$r<r_2$, setting $r':=(r+r_2)/2$ and
$\zeta:=\min\{1,(r'-\rho)^+/(r'-r)\}$,
\[
\int_{\Omega_r}|\nabla_Xu|^2d\tau\le8\|u\|_\infty^2\Big(\int|\nabla_X\zeta|^2
d\tau+1\Big)<\infty .
\]
\end{lemma}

\begin{proof}
For $k$ with $t_k<r$, the function $\psi:=u\,\zeta^2(1-\eta_k)^2$ is
Lipschitz with compact support in $X\cap\{s_k\le\rho\le r'\}\subset
\Omega_{r_2}$, hence $\int\nabla_Xu\cdot\nabla_X\psi\,d\tau=0$. With
$\chi:=\zeta(1-\eta_k)$ we obtain
$\int|\nabla_Xu|^2\chi^2=-2\int u\chi\,\nabla_Xu\cdot\nabla_X\chi\le
\frac12\int|\nabla_Xu|^2\chi^2+2\int u^2|\nabla_X\chi|^2$ and
$|\nabla_X\chi|^2\le2|\nabla_X\zeta|^2+2|\nabla_X\eta_k|^2$, whence
$\int|\nabla_Xu|^2\chi^2\le8\|u\|_\infty^2(\int|\nabla_X\zeta|^2+1/k)$.
Since $\chi=1$ on $A(t_k,r)$ and $t_k\to0$, the claim follows by monotone
convergence.
\end{proof}

\begin{lemma}[Maximum principle]
\label{lemma:massimo}
Let $u$ be harmonic and bounded on $\Omega_{r_2}$. For every $r<r_2$,
\[
\sup_{\Omega_r}u=\max_{L_r}u,\qquad\inf_{\Omega_r}u=\min_{L_r}u .
\]
\end{lemma}

\begin{proof}
Let $M_r:=\max_{L_r}u$, $\varepsilon>0$ and $w:=(u-M_r-\varepsilon)^+$ on
$\Omega_r$. By continuity of $u$ and compactness of $L_r$, $w$ vanishes in
a neighbourhood of $L_r$; moreover $w$ is bounded and, by
Lemma~\ref{lemma:energia}, $\int_{\Omega_r}|\nabla_Xw|^2<\infty$. The
function $(1-\eta_k)w$, extended by $0$, is Lipschitz with compact
support in $\Omega_{r_2}$, hence
\[
0=\int(1-\eta_k)\,\nabla_Xu\cdot\nabla_Xw\,d\tau-\int w\,\nabla_Xu\cdot
\nabla_X\eta_k\,d\tau .
\]
The second integral is bounded by $\|w\|_\infty\|\nabla_Xu\|_{L^2(
\Omega_{t_k})}k^{-1/2}\to0$; the first tends to $\int|\nabla_Xw|^2$ by
dominated convergence. Hence $w$ is constant on the connected set
$\Omega_r$ and, since it vanishes near $L_r$, $w\equiv0$. The same
argument applied to $-u$ gives the statement for the infimum.
\end{proof}

\section{Harnack inequality on annuli and existence of the limit}
\label{sec:harnack}

\begin{definition}
\label{def:harnack}
Let $b\ge2$. $X$ satisfies the \emph{Harnack condition on annuli}
(H$_b$) if there exist $C_H\ge1$ and $r_H\in(0,r_0/b]$ such that, for
every $r\in(0,r_H)$ and every harmonic function $v\ge0$ on $A(r/b,br)$,
\[
\sup_{L_r}v\le C_H\inf_{L_r}v .
\]
Condition (H$_4$) is denoted by (H).
\end{definition}

\begin{theorem}
\label{thm:continuita-H}
Let $X$ have connected link and satisfy (H$_b$) for some $b\ge2$. If $u$
is harmonic and bounded on $\Omega_{r_2}$, $r_1<\min\{r_2,r_H\}$ and
$r_k:=(2b)^{1-k}r_1$, then
$\omega_k:=\sup_{\Omega_{r_k}}u-\inf_{\Omega_{r_k}}u$ satisfies
\[
\omega_{k+1}\le\frac{C_H-1}{C_H+1}\,\omega_k ,
\]
and $u$ has a limit at $0$.
\end{theorem}

\begin{proof}
Let $M_k:=\sup_{\Omega_{r_k}}u$ and $m_k:=\inf_{\Omega_{r_k}}u$. By
Lemma~\ref{lemma:massimo}, $M_{k+1}=\max_{L_{r_{k+1}}}u$ and
$m_{k+1}=\min_{L_{r_{k+1}}}u$. Since $br_{k+1}<r_k$, the functions
$u-m_k$ and $M_k-u$ are harmonic and nonnegative on
$A(r_{k+1}/b,br_{k+1})$, and (H$_b$) gives
\[
M_{k+1}-m_k\le C_H(m_{k+1}-m_k),\qquad M_k-m_{k+1}\le C_H(M_k-M_{k+1}).
\]
Adding, we obtain $\omega_{k+1}+\omega_k\le C_H(\omega_k-\omega_{k+1})$.
The intervals $[m_k,M_k]$ are nested and their length tends to $0$; their
common point is the limit of $u$.
\end{proof}

Theorem~\ref{thm:continuita-H} uses a Harnack constant independent of
$r$. The following variant allows constants depending on the scale. For
$0<r<r_0/4$ we set
\[
C_H(r):=\sup\Big\{\frac{\sup_{L_r}v}{\inf_{L_r}v}:\ v>0\ \text{harmonic on }
A(r/4,4r)\Big\}\in[1,+\infty],
\]
and we use the convention $1/C_H(r):=0$ when $C_H(r)=+\infty$.
Condition (H) is equivalent to $\sup_{0<r<r_H}C_H(r)<\infty$.

\begin{theorem}
\label{thm:continuita-Hsigma}
Let $X$ have connected link and let $r_k\downarrow0$ be radii with
$r_1<r_0/4$ and $r_{k+1}\le r_k/8$. If
\[
\sum_{k\ge1}\frac1{C_H(r_k)}=+\infty ,
\]
then every harmonic function $u$ bounded on $\Omega_{r_0}$ has a limit
at $0$.
\end{theorem}

\begin{proof}
Let $M_k$, $m_k$ and $\omega_k$ be as in the proof of
Theorem~\ref{thm:continuita-H}, relative to the radii $r_k$, with $b=4$,
and let $C_k:=C_H(r_{k+1})$. Since $4r_{k+1}<r_k$, the functions $u-m_k$
and $M_k-u$ are harmonic and nonnegative on $A(r_{k+1}/4,4r_{k+1})$; each
of them is either identically zero or positive. The computation of
Theorem~\ref{thm:continuita-H} gives, for $C_k<\infty$,
\[
\omega_{k+1}\le\frac{C_k-1}{C_k+1}\,\omega_k=\Big(1-\frac2{C_k+1}\Big)\omega_k ,
\]
and $\omega_{k+1}\le\omega_k$ in any case. Since $C_k\ge1$, we have
$2/(C_k+1)\ge1/C_k$, hence
\[
\omega_{K+1}\le\omega_1\prod_{k=1}^{K}\Big(1-\frac1{C_k}\Big)\le
\omega_1\exp\Big(-\sum_{k=1}^{K}\frac1{C_k}\Big).
\]
By assumption the series $\sum_k1/C_H(r_k)$ diverges, and so does
$\sum_k1/C_k$, which differs from it only in the first term; therefore
$\omega_k\to0$. Since $\sup_{\Omega_r}u-\inf_{\Omega_r}u$ is monotone in
$r$, the claim follows as in Theorem~\ref{thm:continuita-H}.
\end{proof}

\begin{remark}
\label{rem:fonti-oscillazione}
Theorem~\ref{thm:continuita-H} adapts to annuli centred at $0$ the
classical argument by which the Harnack inequality implies the decay of
oscillation (Moser \cite{Moser1961}). The passage to unbounded constants
in Theorem~\ref{thm:continuita-Hsigma}, through the equivalence between
$\prod_k(1-x_k)=0$ and $\sum_kx_k=+\infty$ for $x_k\in[0,1)$, is used in
the same way by Korobenko, Rios, Sawyer and Shen
\cite{KorobenkoRiosSawyerShen2015} for infinitely degenerate equations.
In what follows conditions (H) and (H$_b$) are verified for classes of
germs.
\end{remark}

\section{Volume doubling, Poincar\'e inequality and generalized horns}
\label{sec:VP}

This section introduces a condition on intrinsic balls away from $0$
which is invariant under bi-Lipschitz homeomorphisms of the intrinsic
metric, and derives (H$_b$) from it for germs bi-Lipschitz equivalent to a
cone or to a generalized horn. The Harnack inequality on balls under
volume doubling and Poincar\'e is known (Section~\ref{subsec:VP-H}); the
passage to annuli adapts the argument of Grigor'yan and Saloff-Coste for
remote balls \cite{GrigoryanSaloffCoste2005}; the case of cones falls
within the results of Akutagawa, Carron and Mazzeo (Remark~\ref{rem:ACM}).
The new results of the section are Proposition~\ref{prop:corni-VP} and
Theorem~\ref{thm:corni-generalizzati}.

\subsection{Conditions (VP) and (RA\texorpdfstring{$_b$}{b})}
\label{subsec:VP-def}

\begin{definition}
\label{def:VP}
$X$ satisfies condition \emph{(VP)} if there exist $C_V\ge1$, $a_P>0$,
$N_P\ge2$ and $r_V>0$ such that for every $x\in X^*$ with $|x|<r_V$ and
every $s>0$ with $N_Ps\le|x|/2$
\begin{gather}
\label{eq:VP-raddoppio}
\mathcal H^d\big(B_X(x,2s)\big)\le C_V\,\mathcal H^d\big(B_X(x,s)\big),\\
\label{eq:VP-poincare}
\int_{B_X(x,N_Ps)}|\nabla_Xf|^2\,d\tau\ge\frac{a_P}{s^2}\,
\inf_{\xi\in\mathbb R}\int_{B_X(x,s)}(f-\xi)^2\,d\tau
\qquad\text{for }f\in C^1\big(B_X(x,N_Ps)\big).
\end{gather}
Let $b\ge2$. $X$ satisfies \emph{(RA$_b$)} if there exist $\Lambda_A\ge1$
and $r_A>0$ such that, for $r<r_A$, any two points of $L_r$ are the
endpoints of a curve contained in $A(r/b,br)$ of length at most
$\Lambda_Ar$.
\end{definition}

Both conditions persist when $r_V$ and $r_A$ are decreased; we always
assume $r_V\le r_0/2$ and $r_A\le r_0/(2b)$.
Inequality~\eqref{eq:VP-poincare} is the weak Poincar\'e inequality with
dilation constant $N_P$ \cite{SaloffCoste1992}, restricted to the balls
with $N_Ps\le|x|/2$. Condition (VP) is the analogue near the singular
point of volume doubling and Poincar\'e on remote balls in the sense of
Grigor'yan and Saloff-Coste
\cite[Definition~4.1]{GrigoryanSaloffCoste2005}, with $\{0\}$ in place of
the region at infinity. Since $d_X\ge|\cdot|$, the condition $N_Ps\le|x|/2$
gives $B_X(x,N_Ps)\subset X^*$, with compact closure in $X^*$. The class
of test functions in~\eqref{eq:VP-poincare} extends by density to the
Sobolev space $W^{1,2}(B_X(x,N_Ps))$. Condition (RA$_b$) is modelled on
the relatively connected annuli of Grigor'yan and Saloff-Coste
\begin{NoHyper}\cite[Definition~5.1]{GrigoryanSaloffCoste2005}\end{NoHyper}.

\subsection{Harnack inequality on balls and on annuli}
\label{subsec:VP-H}

The passage from (VP) to condition (H$_b$) takes two steps: a Harnack
inequality on intrinsic balls away from $0$
(Lemma~\ref{lemma:harnack-palle}), then a chain of such balls along the
curves of (RA$_b$) (Proposition~\ref{prop:VP-H}), as in
\cite[Lemma~6.3]{GrigoryanSaloffCoste2005}. The first step needs
two known results, because (VP) does not directly fall under the hypotheses
of the classical theorems of Grigor'yan and Saloff-Coste. The Poincar\'e
inequality in (VP) has an arbitrary dilation constant $N_P$, which cannot
be fixed equal to $2$ because it grows under bi-Lipschitz homeomorphisms
(Proposition~\ref{prop:VP-invariante}); Theorem~\ref{thm:BjornBjorn}
brings it back to $2$ inside a ball. The space $X^*$ is not complete;
Theorem~\ref{thm:LierlSaloffCoste} gives the Harnack inequality on an open
set without any completeness assumption.

\begin{theorem}[A.~Bj\"orn and J.~Bj\"orn {\cite[Theorem~5.1]{BjornBjorn2018}}]
\label{thm:BjornBjorn}
Let $B_0$ be a ball in a metric measure space in which doubling and the
weak Poincar\'e inequality with exponent $2$ and dilation constant
$\lambda_{\mathrm D}$ hold for every ball contained in $B_0$, and the
measure satisfies $\tau(B')/\tau(B)\ge C_\dagger(r_{B'}/r_B)^{\varkappa}$
for balls $B'\subset B\subset B_0$, for some constants
$C_\dagger,\varkappa>0$. If any two points of $B_0$ are joined by a curve
of length at most $L_Q$ times their distance, then the Poincar\'e
inequality holds, with exponent $2$ on both sides and dilation constant
$L_Q$, for every ball $B$ with $\tfrac52L_QB\subset B_0$, with a constant
depending only on the previous constants.
\end{theorem}

\begin{theorem}[Lierl and Saloff-Coste {\cite[Theorem~5.2(i) and Section~5]{LierlSaloffCoste2012}}]
\label{thm:LierlSaloffCoste}
Let $\mathcal O$ be an open subset of a space with a symmetric, strongly
local, regular Dirichlet form whose intrinsic distance induces the
topology, and such that $B(x,r)$ is relatively compact whenever
$B(x,2r)\subset\mathcal O$. If doubling and the Poincar\'e inequality with
dilation constant $2$ hold for the balls $B(x,r)$ with
$B(x,2r)\subset\mathcal O$, then the parabolic Harnack inequality holds
on $\mathcal O$, and hence the elliptic one: for every $\delta\in(0,1)$
there exists $C_E\ge1$ such that
$\sup_{B(x,\delta r)}v\le C_E\inf_{B(x,\delta r)}v$ for every harmonic
$v>0$ on $B(x,r)$, whenever $B(x,2r)\subset\mathcal O$.
\end{theorem}

The constant $C_E$ depends only on $\delta$ and on the doubling and
Poincar\'e constants:
in Saloff-Coste's proof by Moser iteration
\cite[Theorem~3.1]{SaloffCoste1992}, the local Sobolev inequalities
\begin{NoHyper}\cite[Theorem~2.1]{SaloffCoste1992}\end{NoHyper} and the mean value estimate
\begin{NoHyper}\cite[Theorem~3.2]{SaloffCoste1992}\end{NoHyper} have constants depending only on
these constants.

\begin{lemma}
\label{lemma:harnack-palle}
If $X$ satisfies (VP), there exists $C_E\ge1$, depending only on $C_V$,
$a_P$, $N_P$, with the following property. Let $x_0\in X^*$ with
$|x_0|<r_V/2$ and $R_E:=|x_0|/(16N_P^2)$. Every harmonic $v>0$ on
$B_X(x_0,R_E/12)$ satisfies
\[
\sup_{B_X(x_0,R_E/24)}v\le C_E\inf_{B_X(x_0,R_E/24)}v .
\]
\end{lemma}

\begin{proof}
Let $B_0:=B_X(x_0,R_E)$.

\emph{Step 1: radii of the balls contained in $B_0$.} Let
$B_X(x,s)\subset B_0$. Then $|x|>|x_0|-R_E=(16N_P^2-1)R_E\ge63R_E$,
in particular $|x|>|x_0|/2$, and
$|x|<|x_0|+R_E<r_V$. Along the curve of
Lemma~\ref{lemma:curve-radiali} joining $x$ to $0$, the function
$d_X(x,\cdot)$ is continuous and takes every value in $[0,|x|)$; if we
had $s>2R_E$, there would be a point $w$ of the curve with
$2R_E<d_X(x,w)<\min\{s,|x|\}$, hence $w\in B_X(x,s)$ and
$d_X(w,x_0)>R_E$, a contradiction. Therefore $s\le2R_E$,
$N_Ps\le|x_0|/(8N_P)\le|x|/2$ and, by (VP), doubling
and~\eqref{eq:VP-poincare} hold for $B_X(x,s)$. Moreover, with
$k_P:=\lceil\log_2N_P\rceil$, doubling at $x$ applies to the radii
$2^js$, $0\le j<k_P$, because $N_P2^js<N_P^2s\le|x_0|/8\le|x|/2$; hence
\begin{equation}\label{eq:confronto-volumi}
\tau\big(B_X(x,N_Ps)\big)\le\tau\big(B_X(x,2^{k_P}s)\big)\le
C_V^{\,k_P}\,\tau\big(B_X(x,s)\big).
\end{equation} Let moreover
$B'=B_X(x',s')\subset B=B_X(x,s)\subset B_0$. Then $B\subset B_X(x',2s)$
and, by the same argument applied to $B'$, $s'\le2s$; iterating doubling
at $x'$ over the radii $t<2s\le4R_E$, for which
$N_Pt\le|x_0|/(4N_P)\le|x'|/2$, we obtain
$\tau(B')/\tau(B)\ge C_V^{-2}(s'/s)^{\varkappa}$ with
$\varkappa=\log_2C_V$.

\emph{Step 2: Dirichlet space.} Consider the form
$\mathcal E(f):=\int|\nabla_Xf|^2\,d\tau$ on $L^2(X^*,\tau)$, the closure
of the form on smooth functions with compact support in $X^*$; it is
symmetric, strongly local and regular, and its intrinsic distance is the
Riemannian distance of $X^*$. For $y,w\in B_0$ every curve from $y$ to
$w$ through $0$ has length at least $|y|+|w|>6R_E>d_X(y,w)$, and the
curves of length less than $2R_E$ starting at $y$ stay in
$B(0,|y|+2R_E)\subset B(0,r_0)$ if $r_V\le r_0/2$, which is not
restrictive; hence on $B_0$ the intrinsic distance coincides with $d_X$.
The ball $B_0$ has compact closure in $X^*$.

\emph{Step 3: reduction to dilation constant $2$.} We apply
Theorem~\ref{thm:BjornBjorn} to the metric space $(X^*,d_X)$ with the
measure $\tau$. On a Riemannian manifold every upper gradient of a
Lipschitz function $f$ is at least $|\nabla_Xf|$ almost everywhere, and
the mean minimizes $\xi\mapsto\int(f-\xi)^2$. Inequality~\eqref{eq:VP-poincare}
extends by density from $C^1$ to the Newtonian space $N^{1,2}$ of the
balls involved, which on a Riemannian manifold coincides with $W^{1,2}$,
and the minimal weak upper gradient of $f\in N^{1,2}$ equals $|\nabla_Xf|$
almost everywhere; hence (VP) gives the Poincar\'e inequality of
Theorem~\ref{thm:BjornBjorn} inside $B_0$, for all pairs formed by a
function and an upper gradient, with $\lambda_{\mathrm D}=N_P$. The distance $d_X$ is a length distance, and
any two points of $B_0$ are joined in $X^*$ by a curve of length at most
$2d_X$; by Theorem~\ref{thm:BjornBjorn} with $L_Q=2$ the Poincar\'e
inequality with dilation constant $2$ holds for the balls $B$ with
$5B\subset B_0$, with a constant depending only on $C_V$, $a_P$, $N_P$.
Theorem~\ref{thm:BjornBjorn} is formulated with averages. For a ball
$B=B_X(x,s)\subset B_0$, dividing~\eqref{eq:VP-poincare} by $\tau(B)$
and using~\eqref{eq:confronto-volumi} gives
\[
\frac1{\tau(B)}\int_B|f-f_B|^2\,d\tau\le\frac{C_V^{\,k_P}}{a_P}\,s^2\,
\frac1{\tau(N_PB)}\int_{N_PB}|\nabla_Xf|^2\,d\tau ,
\]
where $f_B$ is the mean of $f$ on $B$, which minimizes
$\xi\mapsto\int_B(f-\xi)^2$; this is the averaged inequality required by
Theorem~\ref{thm:BjornBjorn}, with constant depending only on $C_V$,
$a_P$, $N_P$. Conversely, the averaged inequality with dilation
constant $2$ given by Theorem~\ref{thm:BjornBjorn} is the form used in
Theorem~\ref{thm:LierlSaloffCoste}.

\emph{Step 4: Harnack on $\mathcal O$.} Let $\mathcal O:=B_X(x_0,R_E/6)$.
If $B_X(x,2t)\subset\mathcal O$, the argument of Step~1 applied to
$\mathcal O$ gives $2t\le R_E/3$, and $d_X(x,x_0)<R_E/6$; hence
$B_X(x,5t)\subset B_X(x_0,R_E/6+5R_E/6)=B_0$. Doubling and the Poincar\'e
inequality with dilation constant $2$ therefore hold for all balls
$B_X(x,t)$ with $B_X(x,2t)\subset\mathcal O$, and
Theorem~\ref{thm:LierlSaloffCoste} applies to $\mathcal O$. Taking
$x=x_0$, $r=R_E/12$ and $\delta=1/2$, so that $B_X(x_0,2r)=\mathcal O$,
gives the claim, with $C_E$ the constant of
Theorem~\ref{thm:LierlSaloffCoste} for $\delta=1/2$.
\end{proof}

\begin{proposition}
\label{prop:VP-H}
If $X$ satisfies (VP) and (RA$_b$), then it satisfies (H$_{2b}$) with
$C_H\le C_E^{\,n_A}$, where $n_A:=\lceil\Lambda_Ab/\varepsilon_E\rceil+1$,
$\varepsilon_E:=1/(384N_P^2)$ and $C_E$ is the constant of
Lemma~\ref{lemma:harnack-palle}. In particular (VP) and (RA$_2$) imply
(H).
\end{proposition}

\begin{proof}
Let $x_0\in X^*$ with $|x_0|<r_V/2$. Lemma~\ref{lemma:harnack-palle}
gives, for every harmonic $v>0$ on the ball $B_X(x_0,|x_0|/(192N_P^2))$,
\[
\sup_{B_X(x_0,\varepsilon_E|x_0|)}v\le C_E\inf_{B_X(x_0,\varepsilon_E|x_0|)}v .
\]
Let $r<\min\{r_V,r_A\}/(2b)$, let $v\ge0$ be harmonic on
$A(r/(2b),2br)$ and $y,y'\in L_r$. The annulus is connected, so by the
strong maximum principle either $v\equiv0$, and there is nothing to
prove, or $v>0$; we assume $v>0$. Let $\sigma_A$ be the curve of (RA$_b$), contained in
$A(r/b,br)$ and of length at most $\Lambda_Ar$, and let
$p_0=y,p_1,\dots,p_{n}=y'$ be points of $\sigma_A$ with
$d_X(p_{i-1},p_i)<\varepsilon_Er/b$; one can take $n\le n_A$. For every
$i$ we have $r/b\le|p_i|\le br<r_V/2$ and, since
$\big||w|-|p_i|\big|\le d_X(w,p_i)$, the ball $B_X(p_i,|p_i|/(192N_P^2))$ is
contained in $A(r/(2b),2br)$. Since
$p_{i}\in B_X(p_{i-1},\varepsilon_E|p_{i-1}|)$, the previous estimate with
$x_0=p_{i-1}$ gives $v(p_i)\le C_Ev(p_{i-1})$ and
$v(p_{i-1})\le C_Ev(p_i)$; iterating, $v(y')\le C_E^{\,n_A}v(y)$.
\end{proof}

\subsection{Bi-Lipschitz invariance}
\label{subsec:VP-invarianza}

\begin{proposition}
\label{prop:VP-invariante}
Let $\Psi$ be a homeomorphism from a neighbourhood of $0$ in $X$ onto a
neighbourhood of $0$ in a set $Y$ as above, $\Psi(0)=0$, bi-Lipschitz for
$d_X$ and $d_Y$ with constant $K_\Psi$. If $Y$ satisfies (VP), then $X$
satisfies (VP), with
\[
N_P'=2N_PK_\Psi^3C_0^Y,\qquad C_V'=K_\Psi^{2d}C_V^{\,m_\Psi},\qquad
a_P'=a_PK_\Psi^{-2d-4},\qquad m_\Psi:=\lceil\log_2(2K_\Psi^2)\rceil ,
\]
where $C_0^Y$ is the constant of~\eqref{eq:distanza-origine} for $Y$, and
with $r_V'=r_V/(K_\Psi C_0^X)$, reduced if necessary so that the balls
involved lie in the domain of $\Psi$.
\end{proposition}

\begin{proof}
$\Psi$ is differentiable almost everywhere on $X^*$, with Jacobian
between $K_\Psi^{-d}$ and $K_\Psi^d$ and $|\nabla_X(F\circ\Psi)|\le
K_\Psi|\nabla_YF|\circ\Psi$ almost everywhere. Let $x\in X^*$ and
$y=\Psi(x)$. For every $t>0$ we have
$B_Y(y,t/K_\Psi)\subset\Psi(B_X(x,t))\subset B_Y(y,K_\Psi t)$. Hence
\[
\mathcal H^d(B_X(x,2s))\le K_\Psi^d\mathcal H^d(B_Y(y,2K_\Psi s))\le
K_\Psi^dC_V^{m_\Psi}\mathcal H^d(B_Y(y,s/K_\Psi))\le
K_\Psi^{2d}C_V^{m_\Psi}\mathcal H^d(B_X(x,s)),
\]
and, for $f\in C^1(B_X(x,N_PK_\Psi^2s))$ and $F:=f\circ\Psi^{-1}$,
\begin{align*}
\inf_\xi\int_{B_X(x,s)}(f-\xi)^2&\le K_\Psi^d\inf_\xi\int_{B_Y(y,K_\Psi
s)}(F-\xi)^2\le\frac{K_\Psi^{d+2}s^2}{a_P}\int_{B_Y(y,N_PK_\Psi
s)}|\nabla_YF|^2\\
&\le\frac{K_\Psi^{2d+4}s^2}{a_P}\int_{B_X(x,N_PK_\Psi^2s)}|\nabla_Xf|^2 ;
\end{align*}
here~\eqref{eq:VP-poincare} for $Y$ extends by density to Lipschitz
functions. The balls of $Y$ involved have radius at most $2K_\Psi^2s$ in
the doubling estimate and $N_PK_\Psi s$ in the Poincar\'e inequality;
their admissibility in $Y$ follows from $|y|\ge d_Y(y,0)/C_0^Y\ge
|x|/(K_\Psi C_0^Y)$ if $N_P's\le|x|/2$. The integral on the right extends
from $B_X(x,N_PK_\Psi^2s)$ to $B_X(x,N_P's)$.
\end{proof}

\begin{lemma}
\label{lemma:RA-trasporto}
Let $X$, $Y$, $\Psi$ be as in Proposition~\ref{prop:VP-invariante}. If
$Y$ satisfies (RA$_b$), then $X$ satisfies (RA$_{b'}$) with
\[
b'=bK_\Psi^2C_0^XC_0^Y,\qquad
\Lambda_A'=K_\Psi^2C_0^X(2+\Lambda_A),
\]
and with $r_A'=r_A/(bK_\Psi C_0^X)$, reduced if necessary so that the
curves involved lie in the image of $\Psi$.
\end{lemma}

\begin{proof}
Let $y,y'\in L^X_r$, $\hat y=\Psi(y)$, $\hat y'=\Psi(y')$. By
\eqref{eq:distanza-origine} and the bi-Lipschitz property of $\Psi$,
$|\hat y|,|\hat y'|\in[r/(K_\Psi C_0^Y),K_\Psi C_0^Xr]$. We join
$\hat y'$ to a point of $L^Y_{|\hat y|}$ by the curve of
Lemma~\ref{lemma:curve-radiali} for $Y$, which lies between the spheres
of radii $|\hat y'|$ and $|\hat y|$ and has length at most
$2K_\Psi C_0^Xr$, and this point to $\hat y$ by the curve of (RA$_b$) in
$A_Y(|\hat y|/b,b|\hat y|)$, of length at most $\Lambda_AK_\Psi C_0^Xr$.
The image under $\Psi^{-1}$ has length at most $\Lambda_A'r$. For a point
$w$ of it we have
\begin{gather*}
|w|\ge\frac{d_X(w,0)}{C_0^X}\ge\frac{d_Y(\Psi(w),0)}{K_\Psi C_0^X}
\ge\frac{|\Psi(w)|}{K_\Psi C_0^X},\\
|w|\le d_X(w,0)\le K_\Psi d_Y(\Psi(w),0)\le K_\Psi C_0^Y|\Psi(w)| ;
\end{gather*}
since $|\Psi(w)|\in[r/(bK_\Psi
C_0^Y),bK_\Psi C_0^Xr]$, we obtain $w\in A(r/b',b'r)$.
\end{proof}

\subsection{Generalized horns}
\label{subsec:corni}

\begin{definition}
\label{def:corno-generalizzato}
Let $L\subset\mathbb R^d\setminus\{0\}$ be a compact connected analytic
submanifold of dimension $d-1$ and let $\gamma\ge1$
be rational. The \emph{generalized horn} with exponent $\gamma$ and
section $L$ is
\[
X_{\gamma,L}:=\{(z^\gamma y,z)\in\mathbb R^d\times\mathbb R:\ y\in L,\
z\ge0\}.
\]
We set $\Pi(z,y):=(z^\gamma y,z)$, $h_L$ the metric induced on $L$,
$a_L:=\max_L|y|$, $D_L:=\mathrm{diam}_{h_L}L$, $\lambda_L>0$ the first
nonzero eigenvalue of the Laplacian of $(L,h_L)$, and $k_L\ge0$ such that
$\mathrm{Ric}_{h_L}\ge-(d-2)k_L^2$.
\end{definition}

Since $L$ is closed in $\mathbb R^d$ and, near each of its points, is the
zero set of finitely many analytic functions, $L$ is semianalytic, hence
subanalytic; as $\gamma$ is rational, $X_{\gamma,L}$ is subanalytic near
$0$. $\Pi$ is an embedding of $(0,\infty)\times L$ onto $X_{\gamma,L}^*$. For
$\gamma=1$ one obtains a cone; for $L$ an ellipsoid, the cusps of
Example~\ref{ex:cuspidi-ellittiche}, with $\gamma=m/2$; for $L$ a sphere,
the rotational cusps of Section~\ref{subsec:cuspidi-rotazione}. The
volume satisfies
$\mathcal H^d(X_{\gamma,L}\cap B(0,r))\approx r^{1+\gamma(d-1)}$.

\begin{lemma}
\label{lemma:prodotto-L}
Let $Z_L:=(\mathbb R\times L,\,d\upsilon^2+h_L)$. Let
$k_L':=\max\{k_L,1\}$, $R_0:=2k_L'D_L$, let $\mathcal B$ be the volume of
the ball of dimension $d$ and curvature $-1$, and let
$c_{\mathrm B}=c_{\mathrm B}(d)\in(0,1)$ be the constant of Lemma~5.1 of
Buser \cite{Buser1982}. For every $p\in Z_L$, $t>0$ and Lipschitz $f$,
\[
\mathrm{vol}\,B_{Z_L}(p,2t)\le C_{Z,V}\,\mathrm{vol}\,B_{Z_L}(p,t),\qquad
\inf_{\xi}\int_{B_{Z_L}(p,t)}(f-\xi)^2\le C_{Z,P}\,t^2\int_{B_{Z_L}(p,2t)}|\nabla f|^2,
\]
with
\[
C_{Z,V}:=\max\Big\{\sup_{0<\varrho\le R_0}\frac{\mathcal B(2\varrho)}
{\mathcal B(\varrho)},\,8\Big\},\qquad
C_{Z,P}:=\max\Big\{4c_{\mathrm B}^{-2(1+R_0)},\,\frac4{\pi^2},\,
\frac1{4D_L^2\lambda_L}\Big\}.
\]
\end{lemma}

\begin{proof}
$Z_L$ is complete and noncompact, with $\mathrm{Ric}\ge-(d-2)k_L^2\ge
-(d-1)k_L'^2$. By \cite[Section~7]{Buser1982}, estimate (4.3.3) and
Lemma~5.1 of \begin{NoHyper}\cite{Buser1982}\end{NoHyper}, which use the curvature only through the
comparisons (4.2)--(4.3), hold on complete noncompact manifolds; we apply
them to the metric $k_L'^2(d\upsilon^2+h_L)$, which satisfies
$\mathrm{Ric}\ge-(d-1)$.

\emph{Radii $t\le2D_L$.} In the rescaled metric the ball has radius
$\varrho=k_L't\le R_0$. Estimate (4.3.3) of \cite{Buser1982} gives
$\mathrm{vol}\,B_{Z_L}(p,2t)\le\frac{\mathcal B(2\varrho)}{\mathcal
B(\varrho)}\mathrm{vol}\,B_{Z_L}(p,t)$; the supremum in $C_{Z,V}$ is
finite because the ratio tends to $2^d$ as $\varrho\to0$. Metric balls
are star-shaped with respect to the centre in the sense of Section~5 of
\begin{NoHyper}\cite{Buser1982}\end{NoHyper}, and Lemma~5.1, with $r=R=\varrho$, gives for the
Neumann isoperimetric constant
\[
\mathfrak h(D)=\inf\Big\{\frac{\mathrm{vol}_{d-1}(\partial A\cap
\mathrm{int}\,D)}{\mathrm{vol}\,A}:\ A\subset D,\ \mathrm{vol}\,A\le
\tfrac12\mathrm{vol}\,D\Big\}
\]
the bound $\mathfrak h\ge c_{\mathrm B}^{1+\varrho}/\varrho\ge
c_{\mathrm B}^{1+R_0}/\varrho$ in the rescaled metric, that is,
$\mathfrak h(B_{Z_L}(p,t))\ge c_{\mathrm B}^{1+R_0}/t$. Cheeger's
inequality, in the form for relatively compact domains recalled in
\begin{NoHyper}\cite[Section~1]{Buser1982}\end{NoHyper}, gives the Poincar\'e inequality on the same
ball with constant $4\mathfrak h^{-2}\le4c_{\mathrm B}^{-2(1+R_0)}t^2$.

\emph{Radii $t\ge2D_L$.} Let $p=(\upsilon_0,\ell)$ and
$I(t'):=(\upsilon_0-t',\upsilon_0+t')$ for $t'>0$. Since
$d_{Z_L}((\upsilon,\ell'),p)^2\le|\upsilon-\upsilon_0|^2+D_L^2$ and
$d_{Z_L}\ge|\upsilon-\upsilon_0|$, we have
$I(t/2)\times L\subset B_{Z_L}(p,t)\subset I(t)\times L$, because
$t^2/4+D_L^2<t^2$. Hence $t\,\mathrm{vol}\,L\le\mathrm{vol}\,B_{Z_L}(p,t)\le
2t\,\mathrm{vol}\,L$, and doubling holds with constant $8$. Moreover
$Q:=I(t)\times L\subset B_{Z_L}(p,\sqrt{t^2+D_L^2})\subset B_{Z_L}(p,2t)$,
and the first nonzero Neumann eigenvalue of the product $Q$ is
$\min\{\pi^2/(2t)^2,\lambda_L\}$; therefore
\begin{align*}
\inf_\xi\int_{B_{Z_L}(p,t)}(f-\xi)^2&\le\int_Q(f-f_Q)^2\le
\max\Big\{\frac{4t^2}{\pi^2},\frac1{\lambda_L}\Big\}\int_Q|\nabla f|^2\\
&\le\max\Big\{\frac4{\pi^2},\frac1{4D_L^2\lambda_L}\Big\}t^2
\int_{B_{Z_L}(p,2t)}|\nabla f|^2 ,
\end{align*}
where $f_Q$ is the mean over $Q$ and we used $t\ge2D_L$.
\end{proof}

\begin{proposition}
\label{prop:corni-VP}
Let $\gamma>1$ and let $z_0\in(0,1]$ with $\delta_L(z_0)\le\tfrac12$,
where $\delta_L(z):=\gamma^2a_L^2z^{2\gamma-2}+\gamma a_Lz^{\gamma-1}$.
With $K_\gamma:=2^{\gamma+1}$, $X_{\gamma,L}$ satisfies
\begin{enumerate}[label=(\roman*)]
\item (VP) with $r_V=z_0/2$, $N_P=128K_\gamma^2$,
$C_V=K_\gamma^{2d}C_{Z,V}^{\lceil\log_2(2K_\gamma^2)\rceil}$ and
$a_P=K_\gamma^{-2d-4}/C_{Z,P}$;
\item (RA$_2$) with $\Lambda_A=1+2D_L$ and $r_A=z_0$.
\end{enumerate}
In particular $X_{\gamma,L}$ satisfies (H).
\end{proposition}

\begin{proof}
\emph{Comparison with the warped product.} For a vector $(\dot z,\dot y)$
tangent to $(0,\infty)\times L$ at $(z,y)$,
\[
\Pi^*g(\dot z,\dot y)=(1+\gamma^2z^{2\gamma-2}|y|^2)\dot z^2
+2\gamma z^{2\gamma-1}\langle y,\dot y\rangle\dot z+z^{2\gamma}|\dot y|^2 .
\]
Setting $g_0:=dz^2+z^{2\gamma}h_L$, from $|y|\le a_L$,
$|\langle y,\dot y\rangle|\le a_L|\dot y|$ and $2z^{2\gamma-1}|\dot
z||\dot y|\le z^{\gamma-1}(\dot z^2+z^{2\gamma}|\dot y|^2)$ we obtain
$|\Pi^*g-g_0|\le\delta_L(z)\,g_0$; hence for $z\le z_0$ we have
$\tfrac12g_0\le \Pi^*g\le\tfrac32g_0$. Moreover $z\le|\Pi(z,y)|\le
z(1+a_L^2z^{2\gamma-2})^{1/2}\le2z$.

\emph{Reduction to a slice.} Let $x=\Pi(z_c,y_c)$ with $|x|<z_0/2$ and
$s>0$ with $N_Ps\le|x|/2$; then $z_c<z_0/2$ and $s\le z_c/N_P$. The balls
used below have radius at most $\rho_*:=2K_\gamma^2s\le z_c/64$. We set
$g_1:=dz^2+z_c^{2\gamma}h_L$ on the slice $F:=\{|z-z_c|\le z_c/2\}$, on
which $z\le3z_c/2<z_0$, so that the comparison with the warped product
holds. On $F$ we have
$(z/z_c)^{2\gamma}\in[2^{-2\gamma},(3/2)^{2\gamma}]$, hence
$K_\gamma^{-2}g_1\le \Pi^*g\le K_\gamma^2g_1$. We denote by $B_1(x,t)$ the
balls of $(\mathbb R\times L,g_1)$, identified with their images under
$\Pi$. For $0<t\le\rho_*$,
\begin{equation}\label{eq:palle-fetta}
B_1(x,t/K_\gamma)\subset B_X(x,t)\subset B_1(x,K_\gamma t).
\end{equation}
Indeed, a curve starting at $x$ of $g_1$-length less than $t/K_\gamma$
satisfies $|z-z_c|<t\le\rho_*$, hence stays in $F$, where its length for
$\Pi^*g$ is less than $t$; a curve starting at $x$ of length less than $t$
in $X$ does not reach $0$, because $t<|x|$, and satisfies
$|z-z_c|<\sqrt2\,t$, because $\Pi^*g\ge\tfrac12dz^2$; it stays in $F$, where
its $g_1$-length is less than $K_\gamma t$. The manifold
$(\mathbb R\times L,g_1)$ is complete and homothetic to $Z_L$ with
factor $z_c^\gamma$.

\emph{(i).} The doubling and Poincar\'e inequalities of
Lemma~\ref{lemma:prodotto-L} are invariant under homotheties and hold
for $g_1$ with the constants $C_{Z,V}$, $C_{Z,P}$, the Poincar\'e
inequality having dilation constant $2$. We repeat the computation in the
proof of Proposition~\ref{prop:VP-invariante} for the identity map, with
constant $K_\gamma$ and Jacobian between $K_\gamma^{-d}$ and $K_\gamma^d$,
and with $(\mathbb R\times L,g_1)$ in place of $Y$. It uses only the
inclusions~\eqref{eq:palle-fetta} for the balls $B_1(x,2K_\gamma s)$,
$B_1(x,s/K_\gamma)$ and $B_X(x,2K_\gamma^2s)$, of radius at most
$\rho_*$, and gives (i), because $2K_\gamma^2\le N_P$.

\emph{(ii).} Let $r<z_0$ and let $p_j=\Pi(z_j,y_j)\in L_r$, $j=1,2$. Then
$z_j\in[r/2,r]$. The curve $t\mapsto\Pi(t,y_1)$, with $t$ between $z_1$
and $z_2$, has length at most $\sqrt{3/2}\,|z_1-z_2|\le r$; the curve
$\Pi(z_2,\sigma_L(\cdot))$, with $\sigma_L$ a minimizing geodesic of
$(L,h_L)$ from $y_1$ to $y_2$, has length at most
$\sqrt{3/2}\,z_2^\gamma D_L\le2D_Lr$, since $r\le1$. Along both curves
$|\Pi(t,y)|\in(t,2t)$ with $t\in[r/2,r]$, because $0\notin L$ and
$a_L^2t^{2\gamma-2}<3$; the curves therefore lie in $A(r/2,2r)$.

The last assertion follows from Proposition~\ref{prop:VP-H}.
\end{proof}

\begin{proposition}[Cones]
\label{prop:coni-VP}
Let $\mathcal C\subset\mathbb R^N$ be the cone over a smooth, compact,
connected submanifold $\Sigma\subset S^{N-1}$ of dimension $d-1$, with
the induced metric $h_\Sigma$, and let
$D_\Sigma:=\mathrm{diam}_{h_\Sigma}\Sigma$. Then $\mathcal C$ satisfies
(VP) for every $r_V>0$, with $N_P=128K_1^2$, $K_1:=4$, and the constants
$C_V$, $a_P$ of Proposition~\ref{prop:corni-VP} with $\gamma=1$ and $L$
replaced by $\Sigma$; it also satisfies (RA$_2$) with
$\Lambda_A=\max\{1,D_\Sigma\}$ and every $r_A>0$. In particular $\mathcal C$
satisfies (H).
\end{proposition}

\begin{proof}
In the coordinates $x=\rho\vartheta$, $\vartheta\in\Sigma$, the metric of
$\mathcal C$ is $g=d\rho^2+\rho^2h_\Sigma$. We repeat the proof of
Proposition~\ref{prop:corni-VP} with $z=\rho$, $\gamma=1$ and $L$
replaced by $(\Sigma,h_\Sigma)$, taking $g_0=g$: the comparison with the
warped product is an equality, $|x|=\rho$, and the condition on $z_0$
plays no role. On the slice $\{|\rho-\rho_c|\le\rho_c/2\}$ the metric
$g_1:=d\rho^2+\rho_c^2h_\Sigma$ satisfies $K_1^{-2}g_1\le g\le K_1^2g_1$,
and the rest of the proof of (i) is unchanged;
Lemma~\ref{lemma:prodotto-L} applies to $(\Sigma,h_\Sigma)$, which is
compact, connected and smooth. For (RA$_2$), two points
$\rho\vartheta_1,\rho\vartheta_2\in L_\rho$ are joined in $L_\rho$ by the
curve $\rho\,\sigma_\Sigma(\cdot)$, with $\sigma_\Sigma$ a minimizing
geodesic of $(\Sigma,h_\Sigma)$, of length at most $D_\Sigma\rho$.
\end{proof}

\begin{corollary}
\label{cor:bilip-cono}
Let $\Psi$ be a homeomorphism from a neighbourhood of $0$ in $X$ onto a
neighbourhood of the vertex of a cone $\mathcal C$ as in
Proposition~\ref{prop:coni-VP}, with $\Psi(0)=0$, bi-Lipschitz for the
intrinsic distances with constant $K_\Psi$. Then $X$ satisfies (VP),
(RA$_{b'}$) and (H$_{2b'}$), with $b'=2K_\Psi^2C_0^XC_0^{\mathcal C}$, and
every bounded harmonic function on $X^*$ near $0$ extends continuously
to $0$.
\end{corollary}

\begin{proof}
By Proposition~\ref{prop:coni-VP}, Proposition~\ref{prop:VP-invariante},
Lemma~\ref{lemma:RA-trasporto} and Proposition~\ref{prop:VP-H};
continuity follows from Theorem~\ref{thm:continuita-H}.
\end{proof}

\begin{remark}
\label{rem:ACM}
Proposition~\ref{prop:coni-VP} falls within the results of Akutagawa,
Carron and Mazzeo on stratified spaces
\cite[Proposition~3.1 and Theorem~A]{AkutagawaCarronMazzeo2014}. A germ
that is only bi-Lipschitz equivalent to a cone does not carry an iterated
edge metric, and Corollary~\ref{cor:bilip-cono} is obtained here from the
bi-Lipschitz invariance of Proposition~\ref{prop:VP-invariante}; for a
spherical link, H\"older continuity under metrics quasi-isometric to
$|x|^{2\nu}|dx|^2$ is proved by Wan \cite{Wan2018}.
\end{remark}

\begin{example}[Elliptic cusps]
\label{ex:cuspidi-ellittiche}
Let $a_1,\dots,a_d>0$, $m>2$ rational, and
\[
X:=\Big\{(x,z)\in\mathbb R^d\times\mathbb R:\ \sum_{i=1}^d\frac{x_i^2}
{a_i^2}=z^m,\ z\ge0\Big\},\qquad d\ge2 .
\]
Then $X=X_{m/2,L_a}$, with
$L_a:=\{y\in\mathbb R^d:\ \sum_iy_i^2/a_i^2=1\}$, and by
Proposition~\ref{prop:corni-VP} $X$ satisfies (VP), (RA$_2$) and (H). The
ellipsoid $L_a$ is strictly convex and has positive sectional curvature
by the Gauss equation, so one can take $k_{L_a}=0$. Elliptic cusps are
not bi-Lipschitz equivalent to a cone: the intrinsic diameter of $L_r$ is
of order $r^{m/2}=o(r)$, whereas the image of $L_r$ under a bi-Lipschitz
homeomorphism $\Psi$ onto a cone meets every generator at points with
$\rho\ge r/K_\Psi$ and has diameter at least a positive constant times
$r/K_\Psi$, the constant depending on the link of the cone.
\end{example}

\begin{theorem}
\label{thm:corni-generalizzati}
Let $X$ be as above, of dimension $d\ge2$, and let $\Psi$ be a
homeomorphism from a neighbourhood of $0$ in $X$ onto a neighbourhood of
$0$ in $X_{\gamma,L}$, $\Psi(0)=0$, bi-Lipschitz for the intrinsic
distances, with $\gamma\ge1$ and $L$ as in
Definition~\ref{def:corno-generalizzato}. Then every bounded harmonic
function on $X^*$ near $0$ extends continuously to $0$.
\end{theorem}

\begin{proof}
The punctured neighbourhoods of $0$ in $X$ and in $X_{\gamma,L}$ are
homeomorphic, hence the link of $X$ is connected. If $\gamma=1$,
$X_{1,L}$ is the cone over the link $\{(y,1)/|(y,1)|:\ y\in L\}$, which is
analytic, compact and connected, and Corollary~\ref{cor:bilip-cono}
concludes. Let $\gamma>1$. By Proposition~\ref{prop:corni-VP},
$X_{\gamma,L}$ satisfies (VP) and (RA$_2$). By
Proposition~\ref{prop:VP-invariante} $X$ satisfies (VP); by
Lemma~\ref{lemma:RA-trasporto} it satisfies (RA$_{b'}$); by
Proposition~\ref{prop:VP-H} it satisfies (H$_{2b'}$), and
Theorem~\ref{thm:continuita-H} concludes.
\end{proof}

\begin{remark}[The case $d=2$]
\label{rem:dim2}
For $d=2$ and connected link the germ is bi-Lipschitz equivalent, for
the intrinsic metric, to a generalized horn
$X_{\gamma,S^1}=\{x^2+y^2=z^{2\gamma},\ z\ge0\}$ with $\gamma\ge1$
(Birbrair \cite[Theorems~8.1 and~8.3]{Birbrair2008}); the exponent
$\gamma$ is rational because subanalytic functions of one variable admit
Puiseux-type expansions.
Theorem~\ref{thm:corni-generalizzati} therefore recovers removability for
$d=2$ (Mosca \cite{Mosca2026}, Grieser \cite{Grieser1999}) using only the
bi-Lipschitz type of the germ and not the conformal structure.
\end{remark}

\section{Positive harmonic functions: asymptotic decomposition and a B\^ocher-type theorem}
\label{sec:picard}

This section deals with the B\^ocher property of
Section~\ref{subsec:problema} for $d\ge3$. The case of a regular point is
B\^ocher's theorem \cite{Bocher1903} (see Axler, Bourdon and Ramey
\cite{AxlerBourdonRamey1992}).

Throughout, $X\subset\mathbb R^N$ is as in Section~\ref{subsec:problema},
with connected link $L$ and $d=\dim X\ge3$, except in
Section~\ref{subsec:AD-H}, where $d\ge2$.

\subsection{The asymptotic decomposition condition}
\label{subsec:decomposizione}

By the product structure (Section~\ref{subsec:prodotto}) we identify
$X\cap B(0,r_0)\setminus\{0\}$ with $L\times(0,r_0)$. We denote by
$\mathcal P$ the set of positive harmonic functions defined on some
punctured neighbourhood $\Omega_r$ of $0$, $0<r\le r_0$.

\begin{definition}[Condition (AD)]
\label{def:DA}
$X$ satisfies the \emph{asymptotic decomposition} condition (AD) if there
exist $r_G\in(0,r_0]$ and a harmonic function $G_X>0$ on $\Omega_{r_G}$
with $G_X\to+\infty$ uniformly in $y\in L$ as $\rho\to0$ (the
\emph{singular profile}) with the following property: for every
$u\in\mathcal P$, defined on $\Omega_r$, there exist
$r'\in(0,\min\{r,r_G\}]$, $q\ge0$ and a bounded harmonic function
$u_{\mathrm b}$ on $\Omega_{r'}$ such that $u=qG_X+u_{\mathrm b}$ on
$\Omega_{r'}$.
\end{definition}

Condition (AD) is a property of the germ of $X$ at $0$. The decomposition
cannot be required on a fixed punctured neighbourhood: a positive
harmonic function may be unbounded near the outer link, as the Poisson
kernel of the unit ball of $\mathbb R^d$ with pole on the unit sphere
shows in the punctured unit ball.

If $X$ satisfies (H$_b$) for some $b\ge2$, then by
Theorem~\ref{thm:continuita-H} the function $u_{\mathrm b}$ extends
continuously to $0$. Condition (AD) expresses the fact that, modulo
bounded functions, the positive harmonic functions near $0$ form a cone
generated by a single function. The profile is in general not
$\rho^{-(d-2)}$: for rotational cusps it is of order
$\rho^{1-m(d-1)/2}$ (Proposition~\ref{prop:cuspide-rotazione}); this is
why (AD) is stated with a general profile $G_X$.

\paragraph{The conical case.} Let $(L,h_L)$ be a closed (compact and
without boundary), connected, smooth Riemannian manifold of dimension
$d-1$, and let $g=d\rho^2+\rho^2h_L$ exactly on $L\times(0,r_0)$. This is
the situation of a set $X\subset\mathbb R^N$ invariant under dilations and
smooth away from $0$: then $L=X\cap S^{N-1}$ with the metric induced by
the unit sphere. We do not require $L$ to be a sphere, nor any curvature
or topological condition beyond connectedness. Let $\phi_j$ ($j\ge0$) be
an orthonormal basis of eigenfunctions of $\Delta_{h_L}$,
$\Delta_{h_L}\phi_j=-\lambda_j\phi_j$, with
$0=\lambda_0<\lambda_1\le\lambda_2\le\cdots$ ($\lambda_0$ is simple
because $L$ is connected). For $u=\sum_jR_j(\rho)\phi_j(y)$ the equation
$\Delta_Xu=0$ is equivalent to
\[
\rho^2R_j''+(d-1)\rho R_j'-\lambda_jR_j=0,
\qquad
\alpha_j^\pm=\frac{-(d-2)\pm\sqrt{(d-2)^2+4\lambda_j}}2,
\]
where $\alpha_j^\pm$ are the roots of the indicial equation
$\alpha^2+(d-2)\alpha-\lambda_j=0$ and $R_j=\rho^{\alpha_j^\pm}$. We have
$\alpha_0^+=0$, $\alpha_0^-=-(d-2)$, and for $j\ge1$, $\alpha_j^+>0$,
$\alpha_j^-<-(d-2)$. We set $\alpha_1:=\alpha_1^+$.

\begin{proposition}
\label{prop:decomposizione}
If $X$ is a cone with metric $d\rho^2+\rho^2h_L$, where $(L,h_L)$ is a
closed, connected, smooth Riemannian manifold of dimension $d-1\ge2$,
then condition (AD) holds with $G_X=\rho^{-(d-2)}$ and, for every
$u\in\mathcal P$,
\begin{equation}\label{eq:decomposizione}
u(y,\rho)=q\,\rho^{-(d-2)}+u_0+O(\rho^{\alpha_1})\qquad(\rho\to0),
\end{equation}
uniformly in $y\in L$.
\end{proposition}

\begin{proof}
Let $\hat u_j(\rho):=\int_Lu(\cdot,\rho)\phi_j\,d\mathrm{vol}_{h_L}$. Then
$\hat u_j$ solves the above Euler equation, hence
$\hat u_j=c_j^+\rho^{\alpha_j^+}+c_j^-\rho^{\alpha_j^-}$ for $j\ge1$ and
$\hat u_0=A+B\rho^{-(d-2)}$. Since $u>0$,
$|\hat u_j(\rho)|\le\|\phi_j\|_\infty\int_Lu(\cdot,\rho)\,d\mathrm{vol}_{h_L}
=\|\phi_j\|_\infty\mathrm{vol}(L)^{1/2}\hat u_0(\rho)\le C_j(1+\rho^{-(d-2)})$.
For $j\ge1$, $\alpha_j^-<-(d-2)$, hence $c_j^-=0$. Let $\Omega_r$ be the
domain of $u$. For $\rho\le r/2$,
$|\hat u_j(\rho)|=|\hat u_j(r/2)|(2\rho/r)^{\alpha_j^+}\le
|\hat u_j(r/2)|(2\rho/r)^{\alpha_1}$. The coefficients
$\hat u_j(r/2)$ of the smooth function $u(\cdot,r/2)$ on $L$ decay
faster than any power of $\lambda_j$, while $\|\phi_j\|_\infty$ grows at
most like a power of $\lambda_j$: by the Sobolev embedding
$H^\kappa(L)\subset C^0(L)$ with $\kappa>(d-1)/2$ we have
$\|\phi_j\|_\infty\le C\|\phi_j\|_{H^\kappa}\le C(1+\lambda_j)^{\kappa/2}$. Hence
$\sum_{j\ge1}\hat u_j\phi_j=O(\rho^{\alpha_1})$ uniformly in $y$. The
constant mode is $(A+B\rho^{-(d-2)})/\sqrt{\mathrm{vol}(L)}$; we set
$q:=B/\sqrt{\mathrm{vol}(L)}$ and $u_0:=A/\sqrt{\mathrm{vol}(L)}$. We
have $q\ge0$, since otherwise $u\to-\infty$ as $\rho\to0$.
\end{proof}

The coefficient $q$ is related to the conserved flux: for harmonic $u$,
$\mathcal J(r):=\int_{L_r}\partial_{\mathbf n}u\,d\mathcal H^{d-1}$, with
$\mathbf n$ the unit normal to $L_r$ in the direction of increasing
$\rho$, is independent of $r$ by the divergence theorem.

\begin{remark}[Comparison with the literature]
\label{rem:LiWan}
The expansion in the proof of Proposition~\ref{prop:decomposizione} is
Kondrat'ev's method \cite{Kondratev1967} for conical points. The case
$L=S^{d-1}$ with metric $(1+\nu)^2g_{S^{d-1}}$, $\nu>-1$, that is, the
metric $|x|^{2\nu}|dx|^2$ on $\mathbb R^d\setminus\{0\}$, is Theorem~1.1
of Jiayu Li and Fangshu Wan \cite{LiWan2019}, proved by the same method.
\end{remark}

\subsection{A B\^ocher-type theorem}
\label{subsec:dicotomia}

In Lemma~\ref{lemma:bocher-b} we assume condition (AD); in
Lemma~\ref{lemma:bocher-c} we assume (AD) and removability, which holds
for instance under (H$_b$) by Theorem~\ref{thm:continuita-H}.

\begin{lemma}[Divergence]
\label{lemma:bocher-b}
Let $u\in\mathcal P$ be unbounded on every punctured neighbourhood of $0$.
Then $u\to+\infty$ as $\rho\to0$, uniformly in $y\in L$.
\end{lemma}

\begin{proof}
Let $u=qG_X+u_{\mathrm b}$ on $\Omega_{r'}$ be as in (AD). If $q=0$, then
$u=u_{\mathrm b}$ is bounded on $\Omega_{r'}$, contrary to the assumption; hence $q>0$. Since
$u_{\mathrm b}$ is bounded and $G_X\to+\infty$ uniformly in $y$,
$u\to+\infty$ uniformly in $y$.
\end{proof}

\begin{lemma}[Cancellation of the singular part]
\label{lemma:bocher-c}
Let $u_1,u_2\in\mathcal P$ tend to $+\infty$ at $0$. There exist
$\xi_1,\xi_2\in\mathbb R$, not both zero, such that $\xi_1u_1+\xi_2u_2$
extends continuously to $0$.
\end{lemma}

\begin{proof}
Let $u_i=q_iG_X+u_{\mathrm b,i}$ on $\Omega_{r_i}$ be as in (AD); since
$u_i$ is unbounded near $0$, $q_i>0$. Set $(\xi_1,\xi_2):=(q_2,-q_1)$.
Then $\xi_1u_1+\xi_2u_2=q_2u_{\mathrm b,1}-q_1u_{\mathrm b,2}$ is harmonic
and bounded on $\Omega_{\min\{r_1,r_2\}}$, hence extends continuously to
$0$ by removability.
\end{proof}

\begin{proposition}[B\^ocher-type theorem]
\label{prop:picard}
If $X$ satisfies (AD) and removability, in particular if it satisfies
(AD) and (H$_b$) for some $b\ge2$, then the B\^ocher property holds at
$0$.
\end{proposition}

\begin{proof}
This is the content of Lemmas~\ref{lemma:bocher-b}
and~\ref{lemma:bocher-c}.
\end{proof}

By Propositions~\ref{prop:decomposizione} and~\ref{prop:coni-VP},
Proposition~\ref{prop:picard} applies to cones over a smooth, closed,
connected link. It is applied below to sets with regular spherical
blow-up (Theorem~\ref{thm:DA-R}) and to rotational cusps
(Proposition~\ref{prop:cuspide-rotazione}).

\subsection{The case of regular spherical blow-up}
\label{subsec:blowup}

Let $\mathrm{Bl}_0X$ be the closure of $\{(|x|,x/|x|):\ x\in X^*\}$ in
$[0,\infty)\times S^{N-1}$, and let
$\Sigma_0:=\mathrm{Bl}_0X\cap(\{0\}\times S^{N-1})$, identified with the
link of the tangent cone of $X$ at $0$.

\begin{definition}[Condition (R)]
\label{def:R}
$X$ has \emph{regular spherical blow-up}, (R) for short, if $\Sigma_0$ is
a smooth closed submanifold of $S^{N-1}$ of dimension $d-1$ and there
exist an integer $\ell\ge1$, $\mu_0>0$ and a map
\[
\Theta:[0,\mu_0)\times\Sigma_0\to[0,\infty)\times S^{N-1},\qquad
\Theta(\mu,y)=(\mu^\ell,\Upsilon(\mu,y)),
\]
with $\Upsilon$ of class $C^\infty$ up to $\mu=0$ and $\Upsilon(0,y)=y$,
which is a homeomorphism onto a neighbourhood of $\{0\}\times\Sigma_0$ in
$\mathrm{Bl}_0X$ and restricts to a diffeomorphism of
$(0,\mu_0)\times\Sigma_0$ onto its image.
\end{definition}

Since $X$ is closed near $0$, the spherical lift
$\{(|x|,x/|x|):\ x\in X^*\}$ is closed in $(0,r_0)\times S^{N-1}$; hence
$\mathrm{Bl}_0X\cap\{\rho>0\}$ coincides, near $\rho=0$, with this lift,
and the restriction of $\Theta$ to $\mu>0$ parametrizes the points of
$X^*$ near $0$. Without closedness the closure in the definition of
$\mathrm{Bl}_0X$ could add points that do not belong to $X$.

\begin{remark}
\label{rem:esempi-R}
The following classes satisfy (R).
\begin{enumerate}[label=(\roman*)]
\item Cones over a smooth closed submanifold of $S^{N-1}$, with $\ell=1$
and $\Theta(\mu,y)=(\mu,y)$.
\item Nondegenerate perturbations of a homogeneous cone: $X=\{F=0\}$ with
$F:\mathbb R^N\to\mathbb R^{N-d}$ analytic, $F=F_k+O(|x|^{k+1})$, $F_k$
homogeneous of degree $k\ge1$ and $dF_k$ of rank $N-d$ at every point of
$\{F_k=0\}\setminus\{0\}$. We have
$F(\rho\vartheta)/\rho^k=F_k(\vartheta)+\rho R(\rho,\vartheta)$ with $R$
analytic. By Euler's identity $dF_k(\vartheta)\vartheta=kF_k(\vartheta)=0$
at the points of $\Sigma_0:=\{F_k=0\}\cap S^{N-1}$, so the restriction of
$dF_k(\vartheta)$ to $T_\vartheta S^{N-1}$ still has rank $N-d$; the
implicit function theorem gives (R) with $\ell=1$. This class contains
the real quadratic forms $Q$ with $\nabla Q\ne0$ on
$\{Q=0\}\setminus\{0\}$, as in Example~\ref{ex:log}, and the complex
analytic germs $\{f=0\}\subset\mathbb C^n$ with $f=f_k+O(|z|^{k+1})$ and
smooth projective tangent cone $\{f_k=0\}$, such as $A_1$, taking
$F=(\mathrm{Re}\,f,\mathrm{Im}\,f)$: the complex differential $df_k$ is
nonzero, hence onto $\mathbb C$, at the points of
$\{f_k=0\}\setminus\{0\}$.
\item Puiseux-type graphs $X=\{(x,|x|^{a/\ell}):\ x\in\mathbb R^d\}$, with
$a>\ell\ge1$ coprime integers. The tangent cone is
$\mathbb R^d\times\{0\}$. With $x=r\vartheta$ and $t:=r^{1/\ell}$ we have
$\rho=t^\ell(1+t^{2(a-\ell)})^{1/2}$, and
$\mu:=\rho^{1/\ell}=t(1+t^{2(a-\ell)})^{1/(2\ell)}$ is an analytic function
of $t$ with derivative $1$ at $0$, hence invertible. The direction
$(\vartheta,t^{a-\ell})/(1+t^{2(a-\ell)})^{1/2}$ is analytic in
$(\mu,\vartheta)$ and equals $(\vartheta,0)$ for $\mu=0$; thus (R) holds
with this $\ell$.
\end{enumerate}
Condition (R) fails for generalized horns with $\gamma>1$, in particular
for elliptic and rotational cusps, whose tangent cone is a half-line, and
for complex singularities with singular projective tangent cone, such as
$A_k$ with $k\ge2$ in $\mathbb C^n$, $n\ge3$. Since a set satisfying (R) is bi-Lipschitz equivalent
to its tangent cone (proof of Theorem~\ref{thm:DA-R}), and a generalized
horn with $\gamma>1$ is not bi-Lipschitz equivalent to a cone, by the
argument of Example~\ref{ex:cuspidi-ellittiche} applied to $L_r$, whose
intrinsic diameter is $O(r^\gamma)$, the classes of
Theorems~\ref{mthm:A} and~\ref{mthm:C} are disjoint.
\end{remark}

Under (R) we have $\rho=\mu^\ell$; we denote by $h_0$ the metric induced
by $S^{N-1}$ on $\Sigma_0$, and $\Sigma_0$ is connected, being
diffeomorphic to the link of $X$. The point of $X$ with coordinates
$(\mu,y)$ is $\rho\,\Upsilon$; since $|\Upsilon|=1$ implies
$\Upsilon\cdot d\Upsilon=0$,
\begin{equation}\label{eq:metrica-b}
\rho^{-2}g=\Big(\frac{d\rho}\rho\Big)^2+|d\Upsilon|^2
=\ell^2\Big(\frac{d\mu}\mu\Big)^2+h_0+\mu\,\mathcal K ,
\end{equation}
where $\mathcal K$ is a combination of $(d\mu/\mu)^2$, $(d\mu/\mu)\,dy^i$
and $dy^idy^j$ with coefficients smooth up to $\mu=0$: here we use
$d\Upsilon=\mu\,\partial_\mu\Upsilon\,(d\mu/\mu)+d_y\Upsilon$ and
$|d_y\Upsilon(0,y)|^2=h_0$. Hence $\rho^{-2}g$ is a b-metric, that is, a
smooth metric on the b-tangent bundle in the sense of Melrose
\cite[Chapter~2]{Melrose1993}, on $[0,\mu_0)\times\Sigma_0$; it is exact
in the sense of \begin{NoHyper}\cite[Section~2.3]{Melrose1993}\end{NoHyper} only up to the term
$\mu\,\mathcal K$, which plays no role below. At the boundary it
coincides
with the metric of the cone over $(\Sigma_0,h_0)$. Let $\alpha_j^\pm$ be the exponents of the conical case
computed on the link $(\Sigma_0,h_0)$.

\begin{theorem}
\label{thm:DA-R}
If $d\ge3$ and $X$ satisfies (R), then $X$ satisfies (H$_b$) for some
$b\ge2$ and (AD). The profile $G_X$ has the expansion
\[
G_X=\rho^{-(d-2)}+\sum_{0<k<\ell(d-2)}\mu^{k-\ell(d-2)}\,c_k(y)
+c_{\log}\log\rho+O(1),
\]
with $c_k\in C^\infty(\Sigma_0)$ and $c_{\log}\in\mathbb R$ determined by
the metric. The term $c_{\log}\log\rho$ arises from the term of order
$\mu^{k-\ell(d-2)}$ with $k=\ell(d-2)$, whose exponent coincides with the
indicial root $0$.
\end{theorem}

\begin{proof}
\emph{The operator.} Since $g=\rho^2g_b$ with $g_b:=\rho^{-2}g$, the
formula for a conformal change gives
$\mathsf P:=\rho^2\Delta_X=\Delta_{g_b}+(d-2)\,g_b(d\log\rho,d\,\cdot\,)$.
The Laplacian of a b-metric is an elliptic b-differential operator and
$d\log\rho=\ell\,d\mu/\mu$ is a smooth section of the b-cotangent bundle,
so $\mathsf P$ is an elliptic b-differential operator with coefficients
smooth up to $\mu=0$. Freezing the coefficients at $\mu=0$,
by~\eqref{eq:metrica-b}, the indicial operator is
$(\rho\partial_\rho)^2+(d-2)\rho\partial_\rho+\Delta_{h_0}$, with
$\rho\partial_\rho=\ell^{-1}\mu\partial_\mu$. If $\phi_j$ is an
eigenfunction of $\Delta_{h_0}$ with eigenvalue $-\lambda_j$, it maps
$\rho^\alpha\phi_j$ to $(\alpha^2+(d-2)\alpha-\lambda_j)\rho^\alpha\phi_j$,
and the indicial roots are the $\alpha_j^\pm$.

\emph{Condition (H$_b$) and growth.} Let
$g_{\mathcal C}:=\rho^2\big(\ell^2(d\mu/\mu)^2+h_0\big)$ be the metric of
the cone $\mathcal C$ over $(\Sigma_0,h_0)$ in the coordinates $(\mu,y)$.
Since the coefficients of $\mathcal K$ are smooth on the compact set
$[0,\mu_1]\times\Sigma_0$, \eqref{eq:metrica-b} gives
$(1-C\mu)g_{\mathcal C}\le g\le(1+C\mu)g_{\mathcal C}$ for $\mu\le\mu_1$.
Consider the map $\Phi_R:\rho\,\Upsilon(\mu,y)\mapsto\rho\,y$, extended by
$\Phi_R(0)=0$; it is a homeomorphism from a neighbourhood of $0$ in $X$
onto a neighbourhood of the vertex of $\mathcal C$. Lengths of curves
avoiding $0$ are multiplied by a factor in
$[(1-C\mu_1)^{1/2},(1+C\mu_1)^{1/2}]$. A curve through $0$ can be
replaced by a curve avoiding $0$ of almost the same length, because the
restriction of $g$ to $L_t$ is comparable to $t^2h_0$ and the intrinsic
diameter of $L_t$ is therefore $O(t)$; the same holds on $\mathcal C$.
Moreover, by Lemma~\ref{lemma:curve-radiali}, a curve of almost minimal
length between two points of $B(0,r)$ stays in $B(0,(1+2C_0)r)$, where
the comparison holds once $r$ is small. Hence $\Phi_R$ is bi-Lipschitz for
the intrinsic distances, and Corollary~\ref{cor:bilip-cono} gives (VP)
and (H$_b$). For $y\in L_{r/2}$ and $y'\in L_r$, we join $y$ to a point
of $L_r$ by the curve of Lemma~\ref{lemma:curve-radiali}, of length at
most $r$, and then to $y'$ within $L_r$ using (H$_b$). Along the first
curve the points satisfy $|x|\ge r/2$, and a chain of at most
$\lceil2/\varepsilon_E\rceil+1$ balls of
Lemma~\ref{lemma:harnack-palle}, as in the proof of
Proposition~\ref{prop:VP-H}, gives $u(y)\le C_\sharp u(y')$ for every
$u\in\mathcal P$, with $C_\sharp$ independent of $r$. Hence
$\sup_{L_{r/2}}u\le C_\sharp\inf_{L_r}u\le C_\sharp\sup_{L_r}u$ and
$u\le C\rho^{-M_\sharp}$ with $M_\sharp=\log_2C_\sharp$.

\emph{Expansion.} We multiply $u\in\mathcal P$ by a cut-off function
$\chi(\mu)$, equal to $1$ near $0$, and extend $\mathsf P$ to a
b-elliptic operator on the compact manifold with boundary
$[0,1]\times\Sigma_0$, joining it to the product b-Laplacian
$(\mu\partial_\mu)^2+\Delta_{h_0}$ by a partition of unity in $\mu$; a
convex combination of second order operators with positive definite
principal symbol is elliptic. The growth estimate gives
$\int_0\mu^{-2s_*}|\chi u|^2\,d\mu/\mu<\infty$ for
$s_*<-\ell M_\sharp$, that is, $\chi u\in\mu^{s_*}L^2_b$, and
$\mathsf P(\chi u)=[\mathsf P,\chi]u$ is smooth and vanishes near the
boundary. By Melrose \cite[Proposition~5.61]{Melrose1993}, $\chi u$ is
polyhomogeneous:
\[
u\sim\sum_{(s,p)}\mu^s(\log\mu)^p\,a_{s,p}(y),\qquad
a_{s,p}\in C^\infty(\Sigma_0),
\]
with exponents $s\in\{\ell\alpha_j^\pm\}+\mathbb N_0$, finitely many below
any threshold, powers of $\log\mu$ given by the index sets
\begin{NoHyper}\cite[(5.118)--(5.121)]{Melrose1993}\end{NoHyper}, and remainders uniform in $y$.

\emph{Leading term.} Let $s_0$ be the smallest exponent with nonzero
coefficient and $p_0$ the highest power of the logarithm with exponent
$s_0$. The coefficient of $\mu^{s_0}(\log\mu)^{p_0}$ in $\mathsf Pu=0$ is
the indicial operator applied to $a_{s_0,p_0}$, hence
$s_0=\ell\alpha_j^\pm$ and $a_{s_0,p_0}$ is an eigenfunction with
eigenvalue $\lambda_j$. If $\lambda_j>0$, then $a_{s_0,p_0}$ has zero
mean and changes sign, so that
$a_{s_0,p_0}(y_1)(\log\mu)^{p_0}<0$ for some $y_1$ and all small $\mu$,
and $u(\mu,y_1)<0$ for small $\mu$, contradicting $u>0$. Hence $a_{s_0,p_0}$ is constant and
$s_0\in\{0,-\ell(d-2)\}$. The coefficient of
$\mu^{s_0}(\log\mu)^{p_0-1}$ is $p_0\,c'\,a_{s_0,p_0}+I_0a_{s_0,p_0-1}$,
where $I_0$ is the indicial operator at the exponent $s_0$, self-adjoint
with kernel the constants, and $c'=2s_0/\ell^2+(d-2)/\ell\ne0$ because
the two roots of the constant mode are simple for $d\ge3$; the range of
$I_0$ is orthogonal to the constants, hence $p_0=0$. It follows that $u$
is bounded, if $s_0=0$, or $u=q\,\rho^{-(d-2)}(1+o(1))$ with $q>0$ and no
exponent smaller than $-\ell(d-2)$.

\emph{The profile.} By Lemma~\ref{lemma:formal} there exists a
polyhomogeneous function
\[
w=\rho^{-(d-2)}+\sum_{0<k<\ell(d-2)}\mu^{k-\ell(d-2)}c_k(y)
+c_{\log}\log\rho+w_+ ,
\]
with $w_+$ bounded, such that $\mathsf Pw$ vanishes to
infinite order at $\mu=0$. Let $f:=-\Delta_X(\chi w)$, which is bounded
and vanishes to infinite order near $0$, and solve $\Delta_Xv=f$ in
$\Omega_{r_0}$ with $v=0$ on $L_{r_0}$ in the Dirichlet space of finite
energy, by the Lax--Milgram lemma and the Poincar\'e inequality along
the curves $y=\mathrm{const}$. The function $v$ is bounded near
$L_{r_0}$ by boundary regularity and belongs to $L^2$, hence to
$\mu^{-\ell d/2}L^2_b$; by Proposition~5.61 of Melrose it is
polyhomogeneous. Since $\Delta_Xv=f$ vanishes to infinite order at
$\mu=0$, the leading exponent $s_0$ of $v$ is an indicial root, as in
the step \emph{Leading term}. A term $\mu^s(\log\mu)^pa(y)$ with $s\ne0$
contributes to the energy
$\int_0\mu^{2s+\ell(d-2)}|\log\mu|^{2p}\,d\mu/\mu$, which diverges if
$2s+\ell(d-2)\le0$; finite energy therefore gives
$s_0>-\ell(d-2)/2$, hence $s_0\ge0$, since no root lies in
$(-\ell(d-2),0)$. Terms $(\log\mu)^p$ with exponent $0$ have finite
energy and are excluded by the argument on the coefficient of
$\mu^{0}(\log\mu)^{p_0-1}$ in the step \emph{Leading term}, which does not
use positivity: the coefficient of the highest power of $\log\mu$ lies in
the kernel of $\Delta_{h_0}$, hence is constant, and $c'=(d-2)/\ell\ne0$
forces $p_0=0$. Thus $v$ is bounded. We set $G_X:=\chi w+v$, adding a constant if necessary so that $G_X>0$ on
$\Omega_{r_0}$.

\emph{Condition (AD).} Let $u\in\mathcal P$ be defined on $\Omega_r$.
By the step \emph{Leading term}, either $s_0=0$, and then $u$ is bounded
on some punctured neighbourhood of $0$ and we take $q=0$ and
$u_{\mathrm b}:=u$, or $s_0=-\ell(d-2)$, with coefficient $q>0$. In the
second case $u_{\mathrm b}:=u-qG_X$ is harmonic on $\Omega_{\min\{r,r_0\}}$
and polyhomogeneous, with no exponents smaller than $-\ell(d-2)$ and with
zero coefficient at $-\ell(d-2)$. Its leading exponent is an indicial
root, hence nonnegative, and the argument in the step \emph{The profile}
excludes the powers of $\log\mu$ at the exponent $0$: $u_{\mathrm b}$ is
bounded on a punctured neighbourhood of $0$.
\end{proof}

\begin{remark}
\label{rem:fonti-DA-R}
The expansion in Theorem~\ref{thm:DA-R} relies on known techniques:
Kondrat'ev's separation of variables \cite{Kondratev1967} for operators
with asymptotically conical coefficients and Melrose's b-calculus
\cite{Melrose1993}. When $\Sigma_0$ is bi-Lipschitz equivalent to a round
sphere, as in Example~\ref{ex:log} and for the Puiseux-type graphs of
Remark~\ref{rem:esempi-R}, the germ is intrinsically bi-Lipschitz
equivalent to a ball of $\mathbb R^d$, $\Delta_X$ becomes in a chart a
uniformly elliptic operator in divergence form, and condition (AD),
without the expansion of the profile, follows from Theorem~4 of Serrin
and Weinberger \cite{SerrinWeinberger1966}.
\end{remark}

\begin{example}[A profile with a logarithmic term]
\label{ex:log}
If $X=f^{-1}(0)\subset\mathbb R^{d+1}$ with $f$ analytic,
$f=f_m+O(|x|^{m+1})$, $f_m$ homogeneous of degree $m$ and
$\nabla f_m\ne0$ on $f_m^{-1}(0)\setminus\{0\}$, condition (R) holds with
$\ell=1$. The function
$F(\rho,\vartheta):=\rho^{-m}f(\rho\vartheta)=f_m(\vartheta)+O(\rho)$ is
analytic up to $\rho=0$; on $\Sigma_0=f_m^{-1}(0)\cap S^d$ the gradient
of $f_m$ is tangent to the sphere by Euler's relation
$\vartheta\cdot\nabla f_m=mf_m$ and nonzero, and the implicit function
theorem gives (R). For instance, for
$X=\{(x,z)\in\mathbb R^3\times\mathbb R:\ |x|^2=z^2+z^3,\ z>0\}\cup\{0\}$,
which is closed near $0$, we have
$d=3$, $\Sigma_0$ is the sphere of radius $1/\sqrt2$ in $S^3$, and a
direct computation on radial functions gives
\[
G_X=\rho^{-1}+\tfrac{\sqrt2}4\log\rho+O(1).
\]
The logarithmic term is that of Theorem~\ref{thm:DA-R}: on radial
functions $\rho^2\Delta_X=a_2(\rho\partial_\rho)^2+a_1\rho\partial_\rho$
with $a_2=1+O(\rho^2)$ and $a_1=1+\tfrac{\sqrt2}4\rho+O(\rho^2)$, hence
$\rho^2\Delta_X(\rho^{-1})=-\tfrac{\sqrt2}4+O(\rho)$, and
$(\rho\partial_\rho)^2+\rho\partial_\rho$ applied to $c_{\log}\log\rho$
gives $c_{\log}=\tfrac{\sqrt2}4$. In particular $u-q\rho^{-(d-2)}$ is
unbounded for $q>0$: condition (AD) holds with the profile $G_X$ and not
with $\rho^{-(d-2)}$. The rotational cusps $|x|^2=z^m$, $m>2$, do not satisfy (R),
because their tangent cone is a half-line.
\end{example}

\subsection{Rotational cusps}
\label{subsec:cuspidi-rotazione}

Let $X=\{(x,z)\in\mathbb R^d\times\mathbb R:\ |x|^2=z^m,\ z\ge0\}$ with
$m>2$ rational, and let $\beta:=m(d-1)/2$. We write
$x=z^{m/2}\vartheta$, $\vartheta\in S^{d-1}$, and
$F_m(z):=\big(1+\tfrac{m^2}4z^{m-2}\big)^{1/2}$; the metric is the warped
product $F_m(z)^2dz^2+z^mg_{S^{d-1}}$, and the level sets of $\rho$ are
the level sets of $z$. We have $X=X_{m/2,S^{d-1}}$ in the notation of
Definition~\ref{def:corno-generalizzato}.

\begin{proposition}
\label{prop:cuspide-rotazione}
If $d\ge3$, the rotational cusp satisfies (AD) with the profile
\[
G_X=\varsigma(z):=\int_z^{z_1}F_m(t)\,t^{-\beta}\,dt
\sim\frac{z^{1-\beta}}{\beta-1}\qquad(z\to0),
\]
where $z_1>0$ is fixed and small,
and every $u\in\mathcal P$ can be written, near $0$, as $u=q\,G_X+u_{\mathrm b}$
with $q\ge0$ and $u_{\mathrm b}$ bounded. Moreover $X$ satisfies (H), and
the B\^ocher property holds at $0$.
\end{proposition}

\begin{proof}
\emph{Radial variable.} The Laplacian is
$\Delta_X=(F_mz^\beta)^{-1}\partial_z\big(z^\beta F_m^{-1}\partial_z\big)+
z^{-m}\Delta_{S^{d-1}}$. In the variable $\varsigma$, with
$d\varsigma=-F_mz^{-\beta}dz$ and $\varsigma\to+\infty$ as $z\to0$
because $\beta>1$, we have
\[
\Delta_X=z^{-2\beta}\big(\partial_\varsigma^2+\mathcal Q\,\Delta_{S^{d-1}}\big),
\qquad \mathcal Q:=z^{2\beta-m}.
\]
Since $1\le F_m\le2$ near $0$ and $\beta>1$, we have
$\varsigma\asymp z^{1-\beta}$ as $z\to0$; in particular
$\varsigma\ge c\,z^{1-\beta}$, that is,
$z\ge(\varsigma/c)^{-1/(\beta-1)}$, and since $2\beta-m>0$ we obtain
$\mathcal Q=z^{2\beta-m}\ge c'\,\varsigma^{-e}$ with
\[
e:=\frac{2\beta-m}{\beta-1}=\frac{2m(d-2)}{m(d-1)-2},\qquad
2-e=\frac{2(m-2)}{m(d-1)-2}>0 .
\]

\emph{Modes.} Let $u\in\mathcal P$ and let $\{Y_j\}$ be an orthonormal
basis of $L^2(S^{d-1})$ of eigenfunctions,
$-\Delta_{S^{d-1}}Y_j=\lambda_jY_j$, with $Y_0$ constant. The rotations
of $S^{d-1}$ are isometries of $X$, so the spherical mean
$\bar u(\varsigma)$ is harmonic and radial: $\bar u''=0$ and
$\bar u=\bar u_0+q\varsigma$, with $q\ge0$ because $\bar u>0$. For
$j\ge1$ the coefficient $u_j(\varsigma):=\int_{S^{d-1}}u\,Y_j$ solves
$u_j''=\lambda_j\mathcal Q\,u_j$. By Proposition~\ref{prop:corni-VP},
with $L=S^{d-1}$ and $\gamma=m/2$, condition (H) holds, so on every level
set $u\le C_H\inf u\le C_H\bar u$, and
\[
|u_j(\varsigma)|\le\|Y_j\|_\infty\,|S^{d-1}|\,C_H\,\bar u(\varsigma)
\le C_j(1+\varsigma).
\]

\emph{Every nonconstant spherical mode is bounded.} The function $w:=u_j^2$ satisfies
$w''=2u_j'^2+2\lambda_j\mathcal Qw\ge0$ and is convex. If it were
unbounded, we would have $w'>0$ and $w\ge A_0>0$ on a half-line. If
$w\ge A\varsigma^p$ on a half-line, then
$w''\ge2\lambda_jcA\varsigma^{p-e}$ and, integrating twice and
restricting the integral to $[\varsigma/2,3\varsigma/4]$, we obtain
$w\ge A'\varsigma^{p+2-e}$ on a half-line. Since $2-e>0$, after finitely
many steps the exponent exceeds $2$, contradicting
$w\le C_j^2(1+\varsigma)^2$. Hence $w$ is bounded; being convex, it is
nonincreasing, and $|u_j(\varsigma)|\le|u_j(\varsigma_0)|$ for
$\varsigma\ge\varsigma_0$.

\emph{Uniform remainder.} The function $u(\varsigma_0,\cdot)$ is smooth
on $S^{d-1}$, so $|u_j(\varsigma_0)|$ decays faster than any power of
$1+\lambda_j$; moreover $\|Y_j\|_\infty\le C(1+\lambda_j)^{\kappa/2}$ for
$\kappa>(d-1)/2$, by the Sobolev embedding, and the series
$\sum_j(1+\lambda_j)^{\kappa/2}|u_j(\varsigma_0)|$ converges, because the
number of eigenvalues smaller than $\lambda$ grows polynomially (Weyl's
law). Since
for every $\varsigma$ the Fourier series of $u(\varsigma,\cdot)$
converges uniformly,
\[
\sup_{\varsigma\ge\varsigma_0}|u-\bar u|\le\sum_{j\ge1}|u_j(\varsigma_0)|\,
\|Y_j\|_\infty<\infty .
\]
We set $u_{\mathrm b}:=\bar u_0+(u-\bar u)$, which is bounded, and
$u=q\,\varsigma+u_{\mathrm b}$. Finally $\varsigma>0$ is harmonic and
tends to $+\infty$ uniformly, and the B\^ocher property follows from
Proposition~\ref{prop:picard}.
\end{proof}

\subsection{The asymptotic decomposition from the Harnack condition}
\label{subsec:AD-H}

In Sections~\ref{subsec:decomposizione}--\ref{subsec:cuspidi-rotazione}
condition (AD) is obtained from an analysis of the specific geometry. In
this subsection we show that (AD) follows from the Harnack inequality on
the links along a sequence of radii; in particular it holds for every
germ satisfying (H$_b$), hence for the whole class of
Theorem~\ref{mthm:A}. The method is that of Pinchover
\cite{Pinchover1994} for Fuchsian type operators: a scale-invariant
Harnack inequality on annuli \begin{NoHyper}\cite[Lemma~6.3]{Pinchover1994}\end{NoHyper}, the
one-dimensionality of the cone of positive solutions vanishing on the
boundary, that is, the Picard principle
\begin{NoHyper}\cite[Theorem~7.1 and Corollary~7.3]{Pinchover1994}\end{NoHyper}, and the
representation of positive solutions near an isolated singular point by
a monotone construction with Dirichlet problems
\begin{NoHyper}\cite[Lemma~8.4]{Pinchover1994}\end{NoHyper}; for the Picard principle see also Nakai
\cite{Nakai1974}. Here the only input is the following condition.

\begin{definition}
\label{def:HS}
$X$ satisfies \emph{(HS)} if there exist $b\ge2$, $C\ge1$ and radii
$r_k\downarrow0$ with $br_k<r_0$ such that, for every $k$ and every
harmonic function $v>0$ on $A(r_k/b,br_k)$,
$\sup_{L_{r_k}}v\le C\inf_{L_{r_k}}v$.
\end{definition}

Condition (H$_b$) implies (HS), and (HS) holds when
$\liminf_{r\to0}C_H(r)<\infty$. Throughout this subsection we assume
(HS). For $0<R\le r_0$ let $\mathcal K_R$ be the set of harmonic
functions $u>0$ on $\Omega_R$ that extend continuously to
$\Omega_R\cup L_R$ with $u=0$ on $L_R$. For $0<r<R$, let $e_r$ be the
harmonic function on $A(r,R)$, continuous on $\overline{A(r,R)}$, with
$e_r=1$ on $L_r$ and $e_r=0$ on $L_R$; it exists because
$\overline{A(r,R)}$ is a compact manifold with smooth boundary
$L_r\cup L_R$, and $0<e_r<1$ on $A(r,R)$ by the strong maximum principle.
For $u$ continuous on $L_r$ we write $m_u(r):=\min_{L_r}u$ and
$M_u(r):=\max_{L_r}u$.

\begin{lemma}
\label{lemma:HS-comparison}
Let $u\in\mathcal K_R$ and let $r=r_k$ with $br_k<R$. Then
$m_u(r)\,e_r\le u\le C\,m_u(r)\,e_r$ on $A(r,R)$.
\end{lemma}

\begin{proof}
The functions $u-m_u(r)e_r$ and $M_u(r)e_r-u$ are harmonic on $A(r,R)$,
continuous on its closure, nonnegative on $L_r$ and zero on $L_R$; by the
maximum principle they are nonnegative. Since $u>0$ is harmonic on
$A(r/b,br)\subset\Omega_R$, (HS) gives $M_u(r)\le C\,m_u(r)$.
\end{proof}

\begin{lemma}
\label{lemma:HS-comparable}
For $u_1,u_2\in\mathcal K_R$ and $x,x_0\in\Omega_R$,
$u_1(x)/u_2(x)\le C^2\,u_1(x_0)/u_2(x_0)$.
\end{lemma}

\begin{proof}
Let $r=r_k$ with $br_k<R$ and $r_k<\min\{|x|,|x_0|\}$. By
Lemma~\ref{lemma:HS-comparison} at the points $x,x_0\in A(r,R)$, where
$e_r>0$, we have $u_1(x)/u_2(x)\le C\,m_{u_1}(r)/m_{u_2}(r)$ and
$u_1(x_0)/u_2(x_0)\ge m_{u_1}(r)/(C\,m_{u_2}(r))$.
\end{proof}

\begin{proposition}[Picard principle]
\label{prop:HS-picard}
For every $R\in(0,r_0]$ there is $G_R\in\mathcal K_R$ such that
$\mathcal K_R=\{tG_R:\ t>0\}$. Moreover $G_R\to+\infty$ uniformly as
$\rho\to0$.
\end{proposition}

\begin{proof}
\emph{Uniqueness.} Fix $x_0\in\Omega_R$ and $u_1,u_2\in\mathcal K_R$ with
$u_1(x_0)=u_2(x_0)=1$, and let $t^*:=\sup\{t\ge0:\ u_1-tu_2\ge0\}$. By
Lemma~\ref{lemma:HS-comparable}, $t^*\ge C^{-2}$, and $t^*\le1$. The
function $w:=u_1-t^*u_2\ge0$ is harmonic on the connected set $\Omega_R$
and vanishes on $L_R$. If $w\not\equiv0$, then $w>0$ by the strong
maximum principle, $w\in\mathcal K_R$, and
Lemma~\ref{lemma:HS-comparable} gives $w\ge C^{-2}w(x_0)u_2=
C^{-2}(1-t^*)u_2$ with $1-t^*=w(x_0)>0$, which contradicts the definition
of $t^*$. Hence $u_1=u_2$.

\emph{Existence.} Let $f_n:=e_{r_n}/e_{r_n}(x_0)$ on $A(r_n,R)$, $r_n<|x_0|$.
Fix $k$ with $br_k<R$ and $r_k<|x_0|$. For $r_n<r_k/b$, the function
$f_n$ is positive and harmonic on $A(r_k/b,br_k)$, and the proof of
Lemma~\ref{lemma:HS-comparison} gives
$m_{f_n}(r_k)e_{r_k}\le f_n\le C\,m_{f_n}(r_k)e_{r_k}$ on $A(r_k,R)$;
evaluating at $x_0$, $m_{f_n}(r_k)\le1/e_{r_k}(x_0)$, hence
$0<f_n\le C\,e_{r_k}/e_{r_k}(x_0)$ on $A(r_k,R)$. By interior elliptic
estimates and a diagonal argument, a subsequence converges locally
uniformly on $\Omega_R$ to a harmonic function $G_R\ge0$ with
$G_R(x_0)=1$ and $G_R\le Ce_{r_k}/e_{r_k}(x_0)$ on $A(r_k,R)$; thus
$G_R$ extends continuously by $0$ to $L_R$, and $G_R>0$ by the strong
maximum principle.

\emph{Growth.} If $G_R$ were bounded, Lemma~\ref{lemma:massimo} on
$\Omega_r$, $r<R$, would give $\sup_{\Omega_r}G_R=\max_{L_r}G_R\to0$ as
$r\to R$, a contradiction. By the maximum principle on
$\overline{A(r',R)}$, $M_{G_R}(r)\le M_{G_R}(r')$ for $r'<r$; hence
$M_{G_R}(r)\to+\infty$ as $r\to0$, and $m_{G_R}(r_k)\ge M_{G_R}(r_k)/C
\to+\infty$. For $r_{k+1}\le r\le r_k$ the minimum principle on
$\overline{A(r_{k+1},r_k)}$ gives $m_{G_R}(r)\ge\min\{m_{G_R}(r_{k+1}),
m_{G_R}(r_k)\}$, so that $\min_{L_r}G_R\to+\infty$.
\end{proof}

\begin{theorem}
\label{thm:AD-H}
Assume (HS). Then:
\begin{enumerate}[label=(\roman*)]
\item every bounded harmonic function on a punctured neighbourhood of $0$
extends continuously to $0$;
\item for every $u\in\mathcal P$, defined on $\Omega_r$, and every
$R'\in(0,r)$, one has $u=q\,G_{R'}+u_{\mathrm b}$ on $\Omega_{R'}$, with
$q\ge0$ and $0\le u_{\mathrm b}\le\max_{L_{R'}}u$; moreover, for
$R'<R\le r_0$, $G_{R'}=c\,G_R+b'$ on $\Omega_{R'}$ with $c>0$ and $b'$
bounded and harmonic;
\item $X$ satisfies (AD), with $G_X:=G_{r_0}$, and the B\^ocher property
holds at $0$.
\end{enumerate}
In particular, (i)--(iii) hold for every germ satisfying (H$_b$) for some
$b\ge2$: the germs of Theorem~\ref{thm:corni-generalizzati}, the germs
satisfying (R), also for $d=2$ (see the proof of
Theorem~\ref{thm:DA-R}), and the elliptic and rotational
cusps (Example~\ref{ex:cuspidi-ellittiche}).
\end{theorem}

\begin{proof}
(i) Passing to a subsequence we may assume $r_{k+1}\le r_k/(2b)$. The
proof of Theorem~\ref{thm:continuita-H}, with the radii $r_k$ in place of
$(2b)^{1-k}r_1$ and $C$ in place of $C_H$, shows that the oscillation of
a bounded harmonic function on $\Omega_{r_k}$ decreases by the factor
$(C-1)/(C+1)$ from $r_k$ to $r_{k+1}$.

(ii) Let $B:=\max_{L_{R'}}u$. For $0<\varepsilon<R'$ let $h_\varepsilon$
be the harmonic function on $A(\varepsilon,R')$, continuous on its
closure, with $h_\varepsilon=u$ on $L_\varepsilon$ and $h_\varepsilon=0$
on $L_{R'}$. By the maximum principle on $\overline{A(\varepsilon,R')}$,
$0\le h_\varepsilon\le u$ and $0\le u-h_\varepsilon\le B$; moreover, for
$\varepsilon'<\varepsilon$, $h_{\varepsilon'}\le u=h_\varepsilon$ on
$L_\varepsilon$, hence $h_{\varepsilon'}\le h_\varepsilon$ on
$A(\varepsilon,R')$. The decreasing limit $h$ is harmonic on
$\Omega_{R'}$, satisfies $0\le h\le h_{\varepsilon_0}$ on
$A(\varepsilon_0,R')$, hence vanishes continuously on $L_{R'}$, and
$0\le u-h\le B$. By the strong maximum principle and
Proposition~\ref{prop:HS-picard}, either $h\equiv0$ or $h=qG_{R'}$ with
$q>0$. Applied to $u=G_R$, which is unbounded on every punctured
neighbourhood of $0$, this gives $G_R=c'G_{R'}+b''$ on $\Omega_{R'}$ with
$c'>0$ and $b''$ bounded, whence the relation between $G_{R'}$ and $G_R$.

(iii) By (ii) with $R=r_0$, every $u\in\mathcal P$ is of the form
$(qc)G_{r_0}+(qb'+u_{\mathrm b})$ on $\Omega_{R'}$, and $G_{r_0}\to+\infty$
uniformly by Proposition~\ref{prop:HS-picard}. This is (AD); together
with (i), Proposition~\ref{prop:picard} gives the B\^ocher property. The
germs listed satisfy (H$_b$) by Theorem~\ref{thm:corni-generalizzati}
and its proof, by the step \emph{Condition (H$_b$) and growth} of the
proof of Theorem~\ref{thm:DA-R} (that step does not use $d\ge3$, so that
(H$_b$) holds under (R) also for $d=2$), and by Example~\ref{ex:cuspidi-ellittiche}
and Proposition~\ref{prop:corni-VP}.
\end{proof}

\begin{remark}
\label{rem:SW-sfera}
Theorem~\ref{thm:AD-H} gives no information on the size of the profile;
Lemma~\ref{lemma:HS-comparison} only shows that $G_R\asymp m_{G_R}(r)\,e_r$
on $A(r,R)$. When $d\ge3$ and $X$ is intrinsically bi-Lipschitz
equivalent to a cone whose link is bi-Lipschitz equivalent to $S^{d-1}$,
the radial extension of a bi-Lipschitz map from the link to $S^{d-1}$
makes $X$ intrinsically bi-Lipschitz equivalent near $0$ to a ball of
$\mathbb R^d$; in such a chart $\Delta_X$ becomes a uniformly elliptic
operator in divergence form with measurable coefficients, and Theorem~4
of Serrin and Weinberger \cite{SerrinWeinberger1966} gives (AD) with a
profile comparable to $\rho^{2-d}$. When $d-1\ne4$, the condition on the
link is topological: for $d-1\le3$ by the classification of smooth
manifolds homeomorphic to spheres, for $d-1\ge5$ by the existence and
uniqueness of Lipschitz structures \cite{Sullivan1979}. The order of
$G_X$ in general is not determined here.
\end{remark}

\section{Examples}
\label{sec:esempi}

This section contains two examples in dimension $d=3$. The first is a
generalized horn with $\gamma>1$ to which
Theorem~\ref{thm:corni-generalizzati} applies and which is not covered by
the known cases; the second is a cone whose link is not a sphere, for
which Proposition~\ref{prop:decomposizione} can be made explicit.

\subsection{A generalized horn with toroidal link}
\label{subsec:corno-toro}

Let $L\subset\mathbb R^3$ be the torus of revolution
\[
L:=\big\{\big((2+\cos\theta_2)\cos\theta_1,\ (2+\cos\theta_2)\sin\theta_1,\ \sin\theta_2\big):\
\theta_1,\theta_2\in\mathbb R/2\pi\mathbb Z\big\},
\]
and let $\gamma>1$ be rational. The generalized horn
$X:=X_{\gamma,L}\subset\mathbb R^3\times\mathbb R$ has dimension $d=3$. The
torus $L$ is semialgebraic, compact, connected and smooth, and does not
contain $0$, because $|y|\ge1$ on $L$; the map $\Pi(z,y)=(z^\gamma y,z)$ is
semialgebraic on $[0,\infty)\times L$ and analytic for $z>0$, hence $X$
is semialgebraic with an isolated singularity at $0$ and link
diffeomorphic to $L$.

\paragraph{Constants.} We have $a_L=\max_L|y|=3$. The intrinsic diameter
satisfies $D_L\le4\pi$: two points are joined by first following an arc
of the parallel $\theta_2=\mathrm{const}$, of length at most $3\pi$, and
then an arc of a meridian, of length at most $\pi$. The Gauss curvature
is $\cos\theta_2/(2+\cos\theta_2)\in[-1,1/3]$, hence
$\mathrm{Ric}_{h_L}\ge-h_L$ and one can take $k_L=1$. The first nonzero
eigenvalue $\lambda_L$ is positive because $L$ is connected; its value
plays no role in the conclusions below. Finally
$\delta_L(z)=9\gamma^2z^{2\gamma-2}+3\gamma z^{\gamma-1}\le\tfrac5{16}$ if
$3\gamma z^{\gamma-1}\le\tfrac14$, and one can take
$z_0:=\min\{1,(12\gamma)^{-1/(\gamma-1)}\}$.

\paragraph{Conclusions.} By Proposition~\ref{prop:corni-VP}, $X$ satisfies
(VP), (RA$_2$) and (H), with constants explicit in terms of $\gamma$,
$D_L$, $k_L$ and $\lambda_L$; by Theorem~\ref{thm:corni-generalizzati}
every bounded harmonic function on a germ bi-Lipschitz equivalent to $X$
extends continuously to $0$. The germ is not covered by the known cases
recalled in Section~\ref{subsec:noti}. The tangent cone of $X$ at $0$ is
the half-line $\{(0,z):\ z\ge0\}$, because $|z^\gamma y|=o(z)$; hence $X$
does not satisfy (R), and it is not bi-Lipschitz equivalent to a cone, by
the argument of Example~\ref{ex:cuspidi-ellittiche}, nor to a ball,
because the link is not a sphere. By Theorem~\ref{thm:AD-H}, $X$ and
every germ bi-Lipschitz equivalent to it satisfy condition (AD), and the
B\^ocher property holds. The order of the singular profile is not
determined here: the proof of Proposition~\ref{prop:cuspide-rotazione},
which gives it for rotational cusps, uses the invariance under the
rotations of $S^{d-1}$, which is missing here.

\subsection{The cone over a Clifford torus}
\label{subsec:toro}

This example illustrates Proposition~\ref{prop:decomposizione} for a link
that is not a sphere.

\paragraph{Setting.} Let
\[
X=\{x\in\mathbb R^4:\ x_1^2+x_2^2=x_3^2+x_4^2\},
\]
a quadric cone of signature $(2,2)$ and dimension $d=3$. The gradient of
$x_1^2+x_2^2-x_3^2-x_4^2$ vanishes only at $0$, so $X\setminus\{0\}$ is
smooth and $0$ is an isolated singularity. The link
\[
L=X\cap S^3=\{x_1^2+x_2^2=x_3^2+x_4^2=\tfrac12\}
\]
is the Clifford torus $S^1(1/\sqrt2)\times S^1(1/\sqrt2)$, which is
connected.

We write $x=\rho\,\mathbf e$ with
\[
\mathbf e=\tfrac1{\sqrt2}\,(\cos\theta_1,\sin\theta_1,\cos\theta_2,\sin\theta_2),
\qquad\theta_1,\theta_2\in\mathbb R/2\pi\mathbb Z,
\]
so that $|\mathbf e|=1$. From $dx=\mathbf e\,d\rho+\rho\,d\mathbf e$,
$\langle \mathbf e,d\mathbf e\rangle=0$ and
$|d\mathbf e|^2=\tfrac12(d\theta_1^2+d\theta_2^2)$ we obtain
\[
g=d\rho^2+\rho^2h_L,\qquad h_L=\tfrac12\,(d\theta_1^2+d\theta_2^2).
\]

\paragraph{Volume and removability.} Since the coefficient of $d\rho^2$
is $1$, we have $|\nabla_X\rho|=1$. The area of the link is
$\mathrm{vol}(L)=(2\pi/\sqrt2)^2=2\pi^2$, hence
$V(\rho)=2\pi^2\rho^{\,2}=2\pi^2\rho^{\,d-1}$. By
Proposition~\ref{prop:stima-capacita}, with
$\Omega=X\cap\{\rho_1<\rho<\rho_2\}$,
\[
\mathrm{cap}_2(L_{\rho_1},L_{\rho_2};\Omega)\le\Big(\int_{\rho_1}^{\rho_2}\frac{dr}{2\pi^2r^2}\Big)^{-1}
=2\pi^2\Big(\frac1{\rho_1}-\frac1{\rho_2}\Big)^{-1}\xrightarrow[\rho_1\to0]{}0 .
\]
Since $X$ is a cone over the smooth, closed, connected link $L$, it
satisfies (H) (Proposition~\ref{prop:coni-VP}); by
Theorem~\ref{thm:continuita-H} every bounded harmonic function on
$X\setminus\{0\}$ extends continuously to $0$.

\paragraph{Spectrum and decomposition.} We have
$\Delta_{h_L}=2(\partial_{\theta_1}^2+\partial_{\theta_2}^2)$, with
eigenfunctions $e^{i(k\theta_1+l\theta_2)}$, $(k,l)\in\mathbb Z^2$, and
eigenvalues
\[
\lambda_{k,l}=2(k^2+l^2).
\]
The first positive eigenvalue is $\lambda_1=2$, of multiplicity $4$
(eigenfunctions $\cos\theta_i$, $\sin\theta_i$). With $d=3$ the indicial
equation $\alpha^2+\alpha-\lambda=0$ has roots
\[
\alpha_{k,l}^\pm=\frac{-1\pm\sqrt{1+8(k^2+l^2)}}2 .
\]
For $(k,l)=(0,0)$: $\alpha^+=0$ and $\alpha^-=-1=-(d-2)$. For
$\lambda_1=2$: $\sqrt{1+8}=3$, hence $\alpha_1=1$ and $\alpha_1^-=-2$. By
Proposition~\ref{prop:decomposizione}, condition (AD) holds with
$G_X=\rho^{-1}$ and every positive harmonic function on
$X\setminus\{0\}$ has the form
\[
u=q\,\rho^{-1}+u_0+O(\rho),\qquad q\ge0,\ \ u_0\in\mathbb R,
\]
uniformly on the link, and the B\^ocher property holds
(Proposition~\ref{prop:picard}).

\appendix
\section{A formal-solutions lemma}
\label{app:formal}

This appendix proves the existence of the function $w$ used in the step
\emph{The profile} of the proof of Theorem~\ref{thm:DA-R}: a
polyhomogeneous function with prescribed leading term $\rho^{-(d-2)}$
and $\mathsf Pw$ vanishing to infinite order at $\mu=0$. The proof is the
recursive construction underlying Lemma~5.44 of Melrose
\cite{Melrose1993}, with explicit control of the exponents and of the
logarithmic terms.

\subsection{Setting}
\label{sec:setting}

We use the notation of Section~\ref{subsec:blowup}. Let $d\ge3$ and let $X$
satisfy (R) with integer $\ell\ge1$; in the coordinates $(\mu,y)\in
[0,\mu_0)\times\Sigma_0$ we have $\rho=\mu^\ell$ and, by~\eqref{eq:metrica-b},
\[
\rho^{-2}g=\ell^2\Big(\frac{d\mu}\mu\Big)^2+h_0+\mu\,\mathcal K ,
\]
where $\mathcal K$ has coefficients smooth up to $\mu=0$. The operator
$\mathsf P=\rho^2\Delta_X$ is a b-differential operator of order $2$ on
$M:=[0,\mu_0)\times\Sigma_0$, elliptic, with coefficients smooth up to
$\mu=0$. We set
\[
D:=\mu\partial_\mu ,\qquad s_*:=-\ell(d-2),\qquad
\mathsf P_0:=\ell^{-2}D^2+(d-2)\ell^{-1}D+\Delta_{h_0},
\]
so that $\mathsf P_0$ is the indicial operator of $\mathsf P$: it is the
operator obtained by freezing the coefficients of $\mathsf P$ at $\mu=0$,
written in the variable $\rho\partial_\rho=\ell^{-1}D$. Since the
coefficients of $\mathsf P-\mathsf P_0$ are smooth and vanish at
$\mu=0$,
\begin{equation}\label{eq:P-split}
\mathsf P=\mathsf P_0+\mu\,\mathsf R ,
\end{equation}
where $\mathsf R$ is a b-differential operator of order $2$ with
coefficients smooth up to $\mu=0$.

Let $(\phi_j)_{j\ge0}$ be an orthonormal basis of $L^2(\Sigma_0,h_0)$ with
$\Delta_{h_0}\phi_j=-\lambda_j\phi_j$, $0=\lambda_0<\lambda_1\le\cdots$,
$\phi_0$ constant ($\Sigma_0$ is connected). On $\mu^\sigma\phi_j$ the
operator $\mathsf P_0$ acts by multiplication by
\[
I_j(\sigma):=\frac{\sigma^2}{\ell^2}+(d-2)\frac{\sigma}{\ell}-\lambda_j ,
\]
whose roots are $\ell\alpha_j^\pm$. For $d\ge3$ the two roots are
distinct for every $j$, hence simple. Moreover
\begin{equation}\label{eq:roots}
\ell\alpha_j^-\le s_*\ \ (\text{with equality only for }j=0),\qquad
\ell\alpha_j^+\ge0\ \ (\text{with equality only for }j=0),
\end{equation}
so that no root lies in the open interval $(s_*,0)$, and the only root
equal to $0$ or to $s_*$ is that of the constant mode.

We say that a function on $M\setminus\{\mu=0\}$ is
\emph{log-smooth of order $(s,p)$} if it has the form
$\mu^s\sum_{i=0}^p(\log\mu)^i\,a_i(\mu,y)$ with $a_i\in C^\infty(M)$.
Since $D\big(\mu^s(\log\mu)^ia\big)=\mu^s\big(s(\log\mu)^ia+i(\log\mu)^{i-1}a
+(\log\mu)^iDa\big)$, every b-differential operator with coefficients
smooth up to $\mu=0$ maps log-smooth functions of order $(s,p)$ to
log-smooth functions of order $(s,p)$.

\subsection{The model equation}

\begin{lemma}
\label{lemma:model}
Let $s\in\mathbb R$, $p\in\mathbb N_0$ and $e_0,\dots,e_p\in
C^\infty(\Sigma_0)$. There exist $b_0,\dots,b_{p+1}\in C^\infty(\Sigma_0)$
such that
\[
\mathsf P_0\Big(\mu^s\sum_{i=0}^{p+1}(\log\mu)^ib_i\Big)
=\mu^s\sum_{i=0}^{p}(\log\mu)^ie_i .
\]
One can take $b_{p+1}=0$ if $s$ is not an indicial root. If $s$ is an
indicial root, $b_{p+1}$ belongs to the finite-dimensional space
$\mathrm{span}\{\phi_j:\ I_j(s)=0\}$; in particular, if $s=0$, then
$b_{p+1}$ is constant.
\end{lemma}

\begin{proof}
Set $t:=\log\mu$, so that $D=\partial_t$, and expand
$e_i=\sum_je_{i,j}\phi_j$; since $e_i$ is smooth, the coefficients
$e_{i,j}$ decay faster than any power of $1+\lambda_j$. For each $j$ we
look for $U_j(t)=\sum_{i}b_{i,j}t^i$ with
\begin{equation}\label{eq:ODE}
I_j(s+\partial_t)\,U_j=E_j(t):=\sum_{i=0}^pe_{i,j}t^i ,
\end{equation}
because $\mathsf P_0\big(e^{st}U_j(t)\phi_j\big)=e^{st}
\big(I_j(s+\partial_t)U_j\big)\phi_j$. The operator
$I_j(s+\partial_t)=I_j(s)+I_j'(s)\partial_t+\ell^{-2}\partial_t^2$ maps
polynomials of degree at most $q$ to polynomials of degree at most $q$.

If $I_j(s)\ne0$ it is invertible on polynomials of degree at most $p$,
with inverse
\[
\sum_{r=0}^{p}(-1)^r I_j(s)^{-1}\big(I_j(s)^{-1}N\big)^r,\qquad
N:=I_j'(s)\partial_t+\ell^{-2}\partial_t^2 ,
\]
since $N$ is nilpotent on these polynomials. The
coefficients $I_j'(s)=2s/\ell^2+(d-2)/\ell$ and $\ell^{-2}$ do not depend
on $j$, and $|I_j(s)|\ge c(1+\lambda_j)$ for all but finitely many $j$;
hence the coefficients $b_{i,j}$ of $U_j$ satisfy
$|b_{i,j}|\le C(1+\lambda_j)^{-1}\max_{i'}|e_{i',j}|$ for large $j$, and
still decay faster than any power of $1+\lambda_j$.

If $I_j(s)=0$, the root is simple, so $I_j'(s)\ne0$, and
$I_j(s+\partial_t)=\partial_t\big(I_j'(s)+\ell^{-2}\partial_t\big)$. The
operator $I_j'(s)+\ell^{-2}\partial_t$ is invertible on polynomials, and
$\partial_t$ maps polynomials of degree at most $p+1$ onto polynomials of
degree at most $p$; hence~\eqref{eq:ODE} has a polynomial solution of
degree at most $p+1$. This happens only for the finitely many $j$ with
$\lambda_j=s^2/\ell^2+(d-2)s/\ell$.

Let $b_i:=\sum_jb_{i,j}\phi_j$. Since $\|\phi_j\|_{C^k}\le
C_k(1+\lambda_j)^{m_k}$ for some $m_k$, by elliptic estimates and the
Sobolev embedding, the series converge in $C^\infty(\Sigma_0)$, and the
$b_i$ solve the equation. The coefficient $b_{p+1}$ receives
contributions only from the indices $j$ with $I_j(s)=0$. For $s=0$ these
are the indices with $\lambda_j=0$, that is, $j=0$ and $b_{p+1}$ is
constant.
\end{proof}

\subsection{Formal solution with prescribed leading term}

\begin{lemma}[Formal solutions]
\label{lemma:formal}
Under the assumptions of Section~\ref{sec:setting} there exist functions
$c_k\in C^\infty(\Sigma_0)$, $0<k<\ell(d-2)$, a constant
$c_{\log}\in\mathbb R$ and a function $w\in C^\infty(M\setminus\{\mu=0\})$
such that:
\begin{enumerate}[label=(\roman*)]
\item $\mathsf Pw$ vanishes to infinite order at $\mu=0$, together with
all its b-derivatives: for every $N$ and every b-differential operator
$Q$, $Q\mathsf Pw=O(\mu^N)$ uniformly in $y$;
\item $w$ is polyhomogeneous, and
\[
w=\rho^{-(d-2)}+\sum_{0<k<\ell(d-2)}\mu^{k+s_*}c_k(y)+c_{\log}\log\rho
+w_+ ,
\]
where $w_+$ is polyhomogeneous with exponents in $\{0\}\cup(0,\infty)$,
no logarithm at the exponent $0$, and finitely many powers of $\log\mu$ at
each positive exponent; in particular $w_+$ is bounded near $\mu=0$.
\end{enumerate}
\end{lemma}

\begin{proof}
\emph{Step 1: the leading term.} Let $u_0:=\mu^{s_*}=\rho^{-(d-2)}$.
Since $I_0(s_*)=0$, we have $\mathsf P_0u_0=0$, and
by~\eqref{eq:P-split} $\mathsf Pu_0=\mu\,\mathsf R(\mu^{s_*})$ is
log-smooth of order $(s_*+1,0)$.

\emph{Step 2: recursion.} We construct, for $n\ge1$, functions
\[
v_n=\mu^{s_*+n}\sum_{i=0}^{q_n}(\log\mu)^i\,b_{n,i}(y),\qquad
b_{n,i}\in C^\infty(\Sigma_0),
\]
such that $u_n:=u_0+v_1+\dots+v_n$ satisfies: $\mathsf Pu_n$ is
log-smooth of order $(s_*+n+1,p_n)$ for some $p_n\in\mathbb N_0$. For
$n=0$ this is Step~1, with $p_0=0$. Assume it holds for $n-1$, and write
$\mathsf Pu_{n-1}=\mu^{s_*+n}\sum_{i\le p_{n-1}}(\log\mu)^ia_i(\mu,y)$
with $a_i$ smooth. By Lemma~\ref{lemma:model} there is $v_n$ of the
stated form with
\[
\mathsf P_0v_n=-\mu^{s_*+n}\sum_{i\le p_{n-1}}(\log\mu)^ia_i(0,y),
\]
with $q_n\le p_{n-1}+1$, and $q_n\le p_{n-1}$ if $s_*+n$ is not an
indicial root. Then
\[
\mathsf Pu_n=\mu^{s_*+n}\sum_i(\log\mu)^i\big(a_i(\mu,y)-a_i(0,y)\big)
+\mu\,\mathsf Rv_n ,
\]
and both terms are log-smooth of order $(s_*+n+1,\max\{p_{n-1},q_n\})$,
because $a_i(\mu,y)-a_i(0,y)=\mu\,\tilde a_i(\mu,y)$ with $\tilde a_i$
smooth. This completes the induction.

\emph{Step 3: exponents and logarithms.} The exponents of the $v_n$ are
$s_*+n$, $n\ge1$. By~\eqref{eq:roots}, for $1\le n<\ell(d-2)$ the number
$s_*+n\in(s_*,0)$ is not an indicial root; starting from $p_0=0$, the
induction gives $q_n=p_n=0$ for these $n$, so that
$v_n=\mu^{s_*+n}c_n(y)$ with $c_n:=b_{n,0}\in C^\infty(\Sigma_0)$. For
$n=\ell(d-2)$ the exponent is $0$, which is a root only for $j=0$; by
Lemma~\ref{lemma:model} with $p=0$,
$v_{\ell(d-2)}=b_{n,0}(y)+b_{n,1}\log\mu$ with $b_{n,1}$ constant, and
the constant part of the kernel can be chosen equal to $0$. Explicitly,
if $a_{0}(0,\cdot)$ has mean value $\bar a$ on $\Sigma_0$, the constant mode of
the equation reads $\big(I_0'(0)+\ell^{-2}\partial_t\big)\partial_tU_0=
-\bar a\,\mathrm{vol}(\Sigma_0)^{1/2}$ with $I_0'(0)=(d-2)/\ell$, whence
$b_{n,1}=-\ell\bar a/(d-2)$. We set
$c_{\log}:=b_{n,1}/\ell$, so that $b_{n,1}\log\mu=c_{\log}\log\rho$. For
$n>\ell(d-2)$ the exponents $s_*+n$ are positive, and the power of the
logarithm can increase by one only when $s_*+n$ coincides with one of the
positive roots $\ell\alpha_j^+$.

\emph{Step 4: asymptotic summation.} Fix smooth vector fields
$V_1,\dots,V_{n_V}$ on $\Sigma_0$ spanning $T_y\Sigma_0$ at every point
$y$; they exist because $\Sigma_0$ is compact. Let $\mathcal M_n$ be the
finite set of monomials $Q=Z_1\cdots Z_k$, $0\le k\le n$, with
$Z_i\in\{\mu\partial_\mu,V_1,\dots,V_{n_V}\}$, and $\mathcal M:=\bigcup_n
\mathcal M_n$, a countable family. Every b-differential operator of order
$k$ with coefficients smooth up to $\mu=0$ is, on $\{\mu\le\mu_0/2\}$, a
finite combination of elements of $\mathcal M_k$ with coefficients that
are smooth and bounded together with their b-derivatives; it therefore
suffices to estimate the elements of $\mathcal M$. Let
$\psi\in C^\infty(\mathbb R)$ with $\psi=1$ on $(-\infty,1/2]$ and
$\psi=0$ on $[1,\infty)$, and set $\psi_n:=\psi(\mu/\varepsilon_n)$.
Since $(\mu\partial_\mu)^i\psi_n$ is bounded independently of
$\varepsilon_n$, each $v_n$ satisfies $|Q(\psi_nv_n)|\le
C_{n,Q}\,\mu^{s_*+n}(1+|\log\mu|)^{q_n}$ for $Q\in\mathcal M$, with
$C_{n,Q}$ independent of $\varepsilon_n\le1$; one can therefore choose
$\varepsilon_n\in(0,1]$, $\varepsilon_n\downarrow0$, such that
\[
\sup_{0<\mu<\varepsilon_n}\ \mu^{-(s_*+n)+1/2}\,
\big|Q\big(\psi_nv_n\big)\big|\le2^{-n}\qquad\text{for all }Q\in\mathcal M_n .
\]
Set
\[
w:=u_0+\sum_{n\ge1}\psi_n\,v_n .
\]
For each $\mu>0$ the sum is finite, and $w\in
C^\infty(M\setminus\{\mu=0\})$.

\emph{Control of the tail.} Fix $N\ge0$ and $Q\in\mathcal M_k$. Then
\[
w-u_N=\sum_{n\le N}(\psi_n-1)v_n+\sum_{N<n<k}\psi_nv_n+\sum_{n\ge
\max\{k,N+1\}}\psi_nv_n .
\]
The first sum vanishes for $\mu<\min_{n\le N}\varepsilon_n/2$. The second
is a finite sum, and each of its terms satisfies $|Q(\psi_nv_n)|\le
C\mu^{s_*+n}(1+|\log\mu|)^{q_n}\le C'\mu^{s_*+N+1/2}$ for $\mu<1$,
because $n\ge N+1$. For the third, $Q\in\mathcal M_n$ and the choice of
$\varepsilon_n$ give $|Q(\psi_nv_n)|\le2^{-n}\mu^{s_*+n-1/2}\le
2^{-n}\mu^{s_*+N+1/2}$ for $\mu<1$. Summing,
\[
Q(w-u_N)=O(\mu^{s_*+N+1/2})\qquad(\mu\to0)
\]
for every fixed $Q\in\mathcal M$, hence for every b-differential
operator $Q$. Thus
$w$ is polyhomogeneous with the expansion $u_0+\sum_nv_n$, and
\[
\mathsf Pw=\mathsf Pu_N+\mathsf P(w-u_N)=O\big(\mu^{s_*+N+1}
(1+|\log\mu|)^{p_N}\big)+O(\mu^{s_*+N+1/2})
\]
with all b-derivatives, because $\mathsf P$ is a b-differential operator.
Since $N$ is arbitrary, (i) follows. Statement (ii) follows from Step~3,
with $w_+$ collecting the terms with $n\ge\ell(d-2)$ other than
$c_{\log}\log\rho$.
\end{proof}

\begin{remark}[Relation with Lemma~5.44 of Melrose]
\label{rem:melrose}
The proof of Lemma~\ref{lemma:formal} given above is self-contained and
does not rely on this remark, which only indicates a second route and is
included for orientation; the correspondence with the notation of
Melrose's index sets is sketched and not used elsewhere.
Lemma~\ref{lemma:formal} can also be derived from
\cite[Lemma~5.44]{Melrose1993}, applied after extending $\mathsf P$ to an
elliptic b-operator on the compact manifold $[0,1]\times\Sigma_0$ as in
the step \emph{Expansion} of the proof of Theorem~\ref{thm:DA-R}, and after
trivializing the b-half-density bundle by
$|d\mu/\mu\ d\mathrm{vol}_{h_0}|^{1/2}$. Let $\chi$ be a cut-off
function equal to $1$ near $\mu=0$ and set $f:=-\mathsf P(\chi\mu^{s_*})$.
By Step~1, $f$ is polyhomogeneous with index set
$E=\{(s_*+1+n,0):\ n\in\mathbb N_0\}$, up to a smooth function supported
away from $\mu=0$. Lemma~5.44 with $\alpha\in(s_*,s_*+1)$ gives a
polyhomogeneous $u$ with $\mathsf Pu-f$ vanishing to infinite order and
index set $F=E\,\overline\cup\,\widehat E^+(\alpha)$, where
$\widehat E^+(\alpha)$ collects the indicial roots with real part larger
than $\alpha$; by~\eqref{eq:roots} these roots are nonnegative. Then
$w:=\chi\mu^{s_*}+u$ satisfies (i), and the exponents of $u$ below $0$
are those of $E$, without logarithms; at the exponent $0$ the extended
union produces at most one power of $\log\mu$, whose coefficient lies in
the kernel of $\Delta_{h_0}$, hence is constant. This gives (ii). The
proof above makes the bookkeeping of the index sets explicit.
\end{remark}

\bibliographystyle{alpha}
\addcontentsline{toc}{section}{\refname}
\bibliography{Harmonic_on_Singular_Subanalytics_refs_v11_en}

@article{Grieser1999,
  author       = {Grieser, Daniel},
  title        = {Local Geometry of Singular Real Analytic Surfaces},
  journal      = {Trans. Amer. Math. Soc.},
  volume       = {355},
  year         = {2003},
  number       = {4},
  pages        = {1559--1577},
  doi          = {10.1090/S0002-9947-02-03168-9},
  note         = {DOI 10.1090/S0002-9947-02-03168-9. Preprint 1999. \url{https://arxiv.org/abs/math/9901062v1}},
  eprint       = {math/9901062v1},
  archivePrefix= {arXiv}
}

@misc{Mosca2026,
  author        = {Mosca, Flavio},
  title         = {A Problem for Isolated Singularities of Surfaces},
  year          = {2026},
  note          = {Original preprint 1989; revised and annotated edition 2026. \url{https://arxiv.org/abs/2608.19452}},
  eprint        = {2608.19452},
  archivePrefix = {arXiv}
}

@article{Bocher1903,
  author  = {B{\^o}cher, Maxime},
  title   = {Singular points of functions which satisfy partial differential equations of the elliptic type},
  journal = {Bull. Amer. Math. Soc.},
  volume  = {9},
  year    = {1903},
  pages   = {455--465}
}

@article{AxlerBourdonRamey1992,
  author  = {Axler, Sheldon and Bourdon, Paul and Ramey, Wade},
  title   = {B\^ocher's theorem},
  journal = {Amer. Math. Monthly},
  volume  = {99},
  year    = {1992},
  pages   = {51--55}
}

@article{Serrin1964b,
  author  = {Serrin, James},
  title   = {Removable singularities of solutions of elliptic equations},
  journal = {Arch. Rational Mech. Anal.},
  volume  = {17},
  year    = {1964},
  pages   = {67--78}
}

@incollection{Hironaka1977,
  author    = {Hironaka, Heisuke},
  title     = {Stratification and flatness},
  booktitle = {Real and Complex Singularities (Proc. Ninth Nordic Summer School, Oslo, 1976)},
  publisher = {Sijthoff \& Noordhoff},
  address   = {Alphen aan den Rijn},
  year      = {1977},
  pages     = {199--265}
}

@article{LiWan2019,
  author  = {Li, Jiayu and Wan, Fangshu},
  title   = {B\^ocher-type theorem on $n$-dimensional manifolds with conical metric},
  journal = {Proc. Amer. Math. Soc.},
  volume  = {147},
  number  = {10},
  year    = {2019},
  pages   = {4527--4538},
  doi     = {10.1090/proc/14554}
}

@misc{Wan2018,
  author        = {Wan, Fangshu},
  title         = {Laplace's equation on $n$-dimensional singular manifolds},
  year          = {2018},
  note          = {\url{https://arxiv.org/abs/1805.04260}},
  eprint        = {1805.04260},
  archivePrefix = {arXiv}
}

@article{SerrinWeinberger1966,
  author  = {Serrin, James and Weinberger, Hans F.},
  title   = {Isolated singularities of solutions of linear elliptic equations},
  journal = {Amer. J. Math.},
  volume  = {88},
  number  = {1},
  year    = {1966},
  pages   = {258--272},
  doi     = {10.2307/2373060},
  note    = {DOI 10.2307/2373060}
}

@article{Kondratev1967,
  author  = {Kondrat'ev, V. A.},
  title   = {Boundary value problems for elliptic equations in domains with conical or angular points},
  journal = {Trudy Moskov. Mat. Obshch.},
  volume  = {16},
  year    = {1967},
  pages   = {209--292},
  note    = {English translation: Trans. Moscow Math. Soc. 1967, 227--313}
}

@article{Moser1961,
  author  = {Moser, J{\"u}rgen},
  title   = {On {H}arnack's theorem for elliptic differential equations},
  journal = {Comm. Pure Appl. Math.},
  volume  = {14},
  number  = {3},
  year    = {1961},
  pages   = {577--591},
  doi     = {10.1002/cpa.3160140329},
  note    = {DOI 10.1002/cpa.3160140329}
}

@book{Melrose1993,
  author    = {Melrose, Richard B.},
  title     = {The {A}tiyah--{P}atodi--{S}inger Index Theorem},
  series    = {Research Notes in Mathematics},
  volume    = {4},
  publisher = {A K Peters},
  address   = {Wellesley, MA},
  year      = {1993},
  isbn      = {1-56881-002-4}
}

@article{Birbrair2008,
  author  = {Birbrair, Lev},
  title   = {Lipschitz geometry of curves and surfaces definable in o-minimal structures},
  journal = {Illinois J. Math.},
  volume  = {52},
  number  = {4},
  year    = {2008},
  pages   = {1325--1353}
}

@article{Buser1982,
  author  = {Buser, Peter},
  title   = {A note on the isoperimetric constant},
  journal = {Ann. Sci. \'Ecole Norm. Sup. (4)},
  volume  = {15},
  number  = {2},
  year    = {1982},
  pages   = {213--230}
}

@article{Grigoryan1991,
  author  = {Grigor'yan, Alexander A.},
  title   = {The heat equation on noncompact {R}iemannian manifolds},
  journal = {Mat. Sb.},
  volume  = {182},
  number  = {1},
  year    = {1991},
  pages   = {55--87},
  note    = {English transl. Math. USSR Sb. 72 (1992), no.~1, 47--77}
}

@article{SaloffCoste1992,
  author  = {Saloff-Coste, Laurent},
  title   = {A note on {P}oincar\'e, {S}obolev, and {H}arnack inequalities},
  journal = {Int. Math. Res. Not. IMRN},
  volume  = {1992},
  year    = {1992},
  number  = {2},
  pages   = {27--38}
}

@incollection{Sullivan1979,
  author    = {Sullivan, Dennis},
  title     = {Hyperbolic geometry and homeomorphisms},
  booktitle = {Geometric Topology (Proc. Georgia Topology Conf., Athens, Ga., 1977)},
  publisher = {Academic Press},
  address   = {New York},
  year      = {1979},
  pages     = {543--555}
}

@article{Valette2024,
  author  = {Valette, Guillaume},
  title   = {On the {L}aplace equation on bounded subanalytic manifolds},
  journal = {Calc. Var. Partial Differential Equations},
  volume  = {63},
  number  = {9},
  year    = {2024},
  pages   = {Paper No.~236},
  doi     = {10.1007/s00526-024-02854-x}
}

@article{Grigoryan1999,
  author  = {Grigor'yan, Alexander},
  title   = {Analytic and geometric background of recurrence and non-explosion of the {B}rownian motion on {R}iemannian manifolds},
  journal = {Bull. Amer. Math. Soc.},
  volume  = {36},
  year    = {1999},
  pages   = {135--249}
}

@article{Bjorn2006,
  author  = {Bj{\"o}rn, Anders},
  title   = {Removable singularities for bounded $p$-harmonic and quasi(super)harmonic functions on metric spaces},
  journal = {Ann. Acad. Sci. Fenn. Math.},
  volume  = {31},
  year    = {2006},
  pages   = {71--95}
}

@article{LierlSaloffCoste2012,
  author  = {Lierl, Janna and Saloff-Coste, Laurent},
  title   = {The {D}irichlet heat kernel in inner uniform domains: local results, compact domains and non-symmetric forms},
  journal = {J. Funct. Anal.},
  volume  = {266},
  number  = {7},
  year    = {2014},
  pages   = {4189--4235}
}

@article{AkutagawaCarronMazzeo2014,
  author  = {Akutagawa, Kazuo and Carron, Gilles and Mazzeo, Rafe},
  title   = {H{\"o}lder regularity of solutions for {S}chr{\"o}dinger operators on stratified spaces},
  journal = {J. Funct. Anal.},
  volume  = {269},
  number  = {3},
  year    = {2015},
  pages   = {815--840}
}

@article{KorobenkoRiosSawyerShen2015,
  author  = {Korobenko, Lyudmila and Rios, Cristian and Sawyer, Eric and Shen, Ruipeng},
  title   = {Local boundedness, maximum principles, and continuity of solutions to infinitely degenerate elliptic equations with rough coefficients},
  journal = {Mem. Amer. Math. Soc.},
  volume  = {269},
  number  = {1311},
  year    = {2021}
}

@article{GrigoryanSaloffCoste2005,
  author  = {Grigor'yan, Alexander and Saloff-Coste, Laurent},
  title   = {Stability results for {H}arnack inequalities},
  journal = {Ann. Inst. Fourier (Grenoble)},
  volume  = {55},
  number  = {3},
  year    = {2005},
  pages   = {825--890}
}

@article{Nakai1974,
  author  = {Nakai, Mitsuru},
  title   = {Martin boundary over an isolated singularity of rotation free density},
  journal = {J. Math. Soc. Japan},
  volume  = {26},
  number  = {3},
  year    = {1974},
  pages   = {483--507}
}

@article{BjornBjorn2018,
  author  = {Bj{\"o}rn, Anders and Bj{\"o}rn, Jana},
  title   = {Local and semilocal {P}oincar{\'e} inequalities on metric spaces},
  journal = {J. Math. Pures Appl.},
  volume  = {119},
  year    = {2018},
  pages   = {158--192}
}

@article{Pinchover1994,
  author  = {Pinchover, Yehuda},
  title   = {On positive {L}iouville theorems and asymptotic behavior of solutions of {F}uchsian type elliptic operators},
  journal = {Ann. Inst. H. Poincar\'e Anal. Non Lin\'eaire},
  volume  = {11},
  number  = {3},
  year    = {1994},
  pages   = {313--341},
}

\end{document}